\documentclass[aap,citesort,MSNbibl,dvips]{arximspdf}
\usepackage{mathbh}

\doi{10.1214/10-AAP744}
\volume{21}
\issue{5}
\pubyear{2011}
\firstpage{1933}
\lastpage{1964}

\makeatletter

\newtheorem{theorem}{Theorem}[section]
\newtheorem{cor}[theorem]{Corollary}
\newtheorem{lem}[theorem]{Lemma}
\newtheorem{prop}[theorem]{Proposition}
\newproclaim{rem}[theorem]{Remark}

\newcommand{\eps}{\varepsilon}
\newcommand{\esssup}{\operatorname{ess}\sup}
\newcommand{\R}{\mathbb R}
\newcommand{\N}{\mathbb N}
\newcommand{\E}{\mathbb E}

\newcommand{\Q}{\mathbb Q}

\newcommand{\supmathitt}{{}^{\mathit t}\hspace*{-0.8pt}}

\makeatother

\begin{document}
\begin{frontmatter}

\title{Numerical simulation of BSDEs with drivers of quadratic growth}
\runtitle{Numerical simulation of quadratic BSDEs}

\begin{aug}
\author[A]{\fnms{Adrien} \snm{Richou}\corref{}\ead[label=e1]{adrien.richou@univ-rennes1.fr}}
\runauthor{A. Richou}
\affiliation{Universit\'{e} Rennes 1}
\address[A]{IRMAR\\
Universit\'{e} de Rennes 1\\
Campus de Beaulieu\\
35042 Rennes Cedex\\
France\\
\printead{e1}} %adresu isvedimo komanda gale!
\end{aug}

% HISTORY:
\received{\smonth{2} \syear{2010}}
\revised{\smonth{9} \syear{2010}}

% ABSTRACT
%
\begin{abstract}
This article deals with the numerical resolution of Markovian backward
stochastic differential equations (BSDEs) with drivers of quadratic
growth with respect to $z$ and bounded terminal conditions. We first
show some bound estimates on the process $Z$ and we specify the Zhang's
path regularity theorem. Then we give a new time discretization scheme
with a nonuniform time net for such BSDEs and we obtain an explicit
convergence rate for this scheme.
\end{abstract}

% KEYWORDS
%

\end{frontmatter}

%s1 ###
\section{Introduction}
Since the early nineties, there has been an increasing interest for
backward stochastic differential equations (BSDEs). These equations
have a wide range of applications in stochastic control, in finance or
in partial differential equation theory. A particular class of BSDE has
been studied for a few years: BSDEs with drivers of quadratic growth
with respect to the variable $z$. This class arises, for example, in
the context of utility optimization problems with exponential utility
functions or alternatively in questions related to risk minimization
for the entropic risk measure (see, e.g., \cite{Hu-Imkeller-Muller-05}).
Many papers deal with existence and uniqueness of solution for such
BSDEs; we refer the reader to \cite
{Kobylanski-00,Lepeltier-SanMartin-98} when the terminal condition is
bounded and \cite{Briand-Hu-06,Briand-Hu-08,Delbaen-Hu-Richou-09} for
the unbounded case. Our concern is rather related to the simulation of
BSDEs and more precisely time discretization of BSDEs coupled with a
forward stochastic differential equation (SDE). Actually, the design of
efficient algorithms which are able to solve BSDEs in any reasonable
dimension has been intensely studied since the first work of
Chevance \cite{Chevance-97} (see, e.g., \cite
{Zhang-04,Bouchard-Touzi-04,Gobet-Lemor-Warin-05}). But in all these
works, the driver of the BSDE is a Lipschitz function with respect to
$z$ and this assumption plays a key role in their proofs. In a recent
paper, Cheridito and Stadje \cite{Cheridito-Stadje-10} studied
approximation of BSDEs by backward stochastic difference equations
which are based on random walks instead of Brownian motions. They
obtain a convergence result when the driver has a subquadratic growth
with respect to $z$ and they give an example where this approximation
does not converge when the driver has a quadratic growth. To the best
of our knowledge, the only work where the time approximation of a BSDE
with a quadratic growth with respect to $z$ is studied is the one of
Imkeller and Reis \cite{Imkeller-dosReis-09}. Notice that, when the
driver has a specific form (roughly speaking, the driver is a sum of a
quadratic term $z \mapsto C\vert z\vert^2$ and a function
that has a linear
growth with respect to $z$), it is possible to get around the problem
by using an exponential transformation method (see \cite
{Imkeller-Reis-Zhang-10}) or by using results on fully coupled
forward--backward differential equations (see \cite{Delarue-Menozzi-06}).

To explain the ideas of this paper, let us introduce $(X,Y,Z)$ the
solution to the forward--backward system
\begin{eqnarray*}
X_t &=& x +\int_0^t b(s,X_s)\,ds+\int_0^t \sigma(s) \,dW_s,\\
Y_t &=& g(X_T) +\int_t^T f(s,X_s,Y_s,Z_s)\,ds-\int_t^T Z_s \,dW_s,
\end{eqnarray*}
where $g$ is bounded, $f$ is locally Lipschitz and has a quadratic
growth with respect to $z$. A well-known result is that when $g$ is a
Lipschitz function with Lipschitz constant $K_g$, then the process $Z$
is bounded by $C(K_g+1)$ (see Theorem \ref{borne globale sur Z}). So,
in this case, the driver of the BSDE is a Lipschitz function with
respect to $z$ and we are able to use standard results about
discretization of BSDEs. Because of the above observation, this paper
will focus on the case that the terminal function $g$ is not Lipschitz.
To obtain our main results, we will assume that $g$ is an $\alpha$-H\"
{o}lder function but it is also possible to adapt our methods when $g$
is not $\alpha$-H\"{o}lder; for example, Remark \ref
{exemple_indicatrice} deals with the case of an indicator function of a
smooth domain. Let us notice that the time approximation of BSDEs with
an irregular terminal function has already been studied by Gobet and
Makhlouf \cite{Gobet-Makhlouf-09} when the generator is a Lipschitz
function with respect to $z$.

In light of previous observation, a simple idea is to do an
approximation of $(Y,Z)$ by the solution $(Y^N,Z^N)$ to the BSDE
\[
Y^N_t = g_N(X_T) +\int_t^T f(s,X_s,Y^N_s,Z^N_s)\,ds-\int_t^T Z^N_s \,dW_s,
\]
where $g_N$ is a Lipschitz approximation of $g$. Thanks to bounded mean
oscillation martingale (BMO martingale in the sequel) tools, we have an
error estimate for this approximation (see, e.g., \cite
{Imkeller-dosReis-09,Briand-Confortola-08} or Proposition \ref{erreur
premiere approximation EDSR}). For example, if $g$ is $\alpha$-H\"
{o}lder, we are able to obtain the error bound $CK_{g_N}^{
{-\alpha
}/({1-\alpha})}$ (see Proposition \ref{erreur approximation g alpha
holder}). Moreover, we can have an error estimate for the time
discretization of the approximated BSDE thanks to any numerical scheme
for BSDEs with Lipschitz driver. But this error estimate depends on
$K_{g_N}$; roughly speaking, this error is $Ce^{CK_{g_N}^2}n^{-1}$ with
$n$ the number of discretization times. The exponential term results
from the use of Gronwall's inequality. Finally, when $g$ is $\alpha
$-H\"
{o}lder and $K_{g_N}=N$, the global error bound is
%
%e1 ###
%
\begin{equation}
\label{introerreurnaif1}
C\biggl(\frac{1}{N^{{\alpha}/({1-\alpha})}}+\frac
{e^{CN^2}}{n}\biggr).
\end{equation}
So, when $N$ increases, $n^{-1}$ will have to become small very quickly
and the speed of convergence turns out to be bad; if we take $N=
(\frac{C}{\eps}\log n)^{1/2}$ with $0<\eps<1$, then the global
error bound becomes $C_{\eps}(\log n)^{{-\alpha
}/({2(1-\alpha)})}$. The same drawback appears in the work of Imkeller and
Reis \cite{Imkeller-dosReis-09}. Indeed, their idea is to do an
approximation of $(Y,Z)$ by the solution $(Y^N,Z^N)$ to the truncated BSDE
\[
Y^N_t = g(X_T) +\int_t^T f(s,X_s,Y^N_s,h_N(Z^N_s))\,ds-\int_t^T Z^N_s \,dW_s,
\]
where $h_N\dvtx\R^{1\times d}\rightarrow\R^{1\times d}$ is a smooth
modification of the projection on the open Euclidean ball of radius $N$
about $0$. Thanks to several statements concerning the path regularity
and stochastic smoothness of the solution processes, the authors show
that for any $\beta\geq1$, the approximation error is lower than
$C_{\beta}N^{-\beta}$. So they obtain the global error bound
%
%e2 ###
%
\begin{equation}
\label{introerreurnaif2}
C_{\beta}\biggl(\frac{1}{N^{\beta}}+\frac{e^{CN^2}}{n}\biggr)
\end{equation}
and, consequently, the speed of convergence also turns out to be bad;
if we take $N=(\frac{C}{\eps}\log n)^{1/2}$ with $0<\eps<1$,
then the global error bound becomes $C_{\beta,\eps}(\log
n
)^{-\beta/2}$.

Another idea is to use an estimate of $Z$ that does not depend on
$K_g$. So we extend a result of \cite{Delbaen-Hu-Bao-09} which shows
%
%e3 ###
%
\begin{equation}
\label{majorationZintroduction}
\vert Z_t\vert\leq M_1+\frac{M_2}{(T-t)^{1/2}},\qquad
0\leq t <T.
\end{equation}
Let us notice that this type of estimation is well known in the case of
drivers with linear growth as a consequence of the Bismut--Elworthy
formula (see, e.g., \cite{Fuhrman-Tessitore-02}). But in our case, we do
not need to suppose that $\sigma$ is invertible. Then, thanks to this
estimation, we know that when $t<T$, $f(t,\cdot,\cdot,\cdot)$ is a Lipschitz
function with respect to $z$ and the Lipschitz constant depends on $t$.
So we are able to modify the classical uniform time net to obtain a
convergence speed for a modified time discretization scheme for our
BSDE; the idea is to put more discretization points near the final time
$T$ than near $0$. Roughly speaking, our discretization grid is equal to
\[
t_k := T \biggl( 1-\biggl( \frac{\eps}{T} \biggr)^{k/n} \biggr),\qquad
0\leq k \leq n,
\]
with $\eps$ a parameter. But due to technical reasons we need to apply
this modified time discretization scheme to the approximated BSDE
\[
Y_t^{N,\eps} = g_N(X_T)+\int_t^T f^{\eps}(s,X_s,Y_s^{N,\eps
},Z_s^{N,\eps
})\,ds-\int_t^T Z_s^{N,\eps} \,dW_s
\]
with
\[
f^{\eps}(s,x,y,z) := \mathbh{1}_{s\leq T-\eps} f(s,x,y,z)+ \mathbh
{1}_{s > T-\eps} f(s,x,y,0).
\]
Thanks to the estimate (\ref{majorationZintroduction}), we obtain a
speed convergence for the time discretization scheme of this\vadjust{\goodbreak}
approximated BSDE (see Theorem \ref{erreur e3 sur [0,T-eps]}).
Moreover, BMO tools give us again an estimate of the approximation
error (see Proposition \ref{erreur premiere approximation EDSR}).
Finally, if we suppose that $g$ is $\alpha$-H\"{o}lder, we prove that
we can choose properly $N$ and $\eps$ to obtain the global error
estimate $Cn^{-{2\alpha}/({(2-\alpha)(2+K)-2+2\alpha})}$ (see
Theorem \ref{erreur_globale}) where $K>0$ depends on constant $M_2$
defined in equation (\ref{majorationZintroduction}) and constants
related to $f$. Let us notice that such a speed of convergence where
constants related to $f$, $g$, $b$ and $\sigma$ appear in the power of
$n$ is unusual. Even if we have an error far better than (\ref
{introerreurnaif1}) or (\ref{introerreurnaif2}), this result is not
very interesting in practice because the speed of convergence strongly
depends on $K$. But, when $b$ is bounded, we prove that we can take
$M_2$ as small as we want in (\ref{majorationZintroduction}). Finally,
we obtain a global error estimate lower than $C_{\eta}n^{-(\alpha
-\eta
)}$ for all $\eta>0$ (see Theorem \ref{vitesse_optimale}).

To conclude, it could be interesting to do some comparisons between our
work and the article of Gobet and Makhlouf \cite{Gobet-Makhlouf-09}. We
already explain that this paper studies the time approximation of
Lipschitz BSDEs with irregular terminal functions. These authors show
that the error of approximation is lower than $C_{\eta}n^{-\alpha}$
when $g$ is an $\alpha$-H\"{o}lder function and the discretization grid
is uniform. So, our better speed of convergence is very close to their
result. Nevertheless, they also show that it is possible to obtain the
classical speed of convergence, that is to say $Cn^{-1}$, when we use
the nonuniform grid given by
\[
t_k:= T-T\biggl(1-\frac{k}{n}\biggr)^{1/\beta},\qquad 0 \leq k \leq n,
\]
with $\beta< \alpha$. It is interesting to notice that we both use
nonuniform time discretization points but their grid is different than
our grid; the accumulation speed of discretization points near the
terminal time $T$ is not the same; it is faster in our case.

The paper is organized as follows. In the introductory Section \ref{sec2} we
recall some of the well-known results concerning SDEs and BSDEs. In
Section \ref{sec3} we establish some estimates concerning the process $Z$; we
show a first uniform bound for $Z$, then a time dependent bound and
finally we specify the classical path regularity theorem. In Section
\ref{sec4}
we define a modified time discretization scheme for BSDEs with a
nonuniform time net and we obtain an explicit error bound.

%s2 ###
\section{Preliminaries}\label{sec2}
%s2.1 ###
\subsection{Notation}\label{sec2.1}
Throughout this paper, $(W_t)_{t \geq0}$ will denote a $d$-dimensional
Brownian motion, defined on a probability space $(\Omega
,\mathcal
{F}, {\mathbb P})$. For $t \geq0$, let $\mathcal{F}_t$ denote the
$\sigma
$-algebra $\sigma(W_s; 0\leq s\leq t)$, augmented with the ${\mathbb P}
$-null sets of $\mathcal{F}$. The Euclidean norm on $\R^d$ will be
denoted by \mbox{$|\cdot|$}. The operator norm induced by \mbox{$|\cdot|$} on the
space of
linear operator is also denoted by \mbox{$|\cdot|$}. %For all probability
%measure $
%expectation with respect to
For $p \geq2$, $m \in\N$, we denote further:
\begin{longlist}[(1)]
\item[(1)] $\mathcal{S}^p(\R^m)$ or $\mathcal{S}^p$ when no confusion is
possible, the space of all adapted processes $(Y_t)_{t \in[0,T]}$ with
values in $\R^m$ normed by
\[
\Vert Y\Vert_{\mathcal{S}^p}=\E\Bigl[\Bigl(\sup_{t \in[0,T]}
\vert Y_t\vert\Bigr)^p\Bigr]^{1/p};
\]
$\mathcal{S}^{\infty}(\R^m)$ or $\mathcal{S}^{\infty}$, the space of
bounded measurable processes;
\item[(2)] $\mathcal{M}^p(\R^m)$ or $\mathcal{M}^p$, the space of all
progressively measurable processes $(Z_t)_{t \in[0,T]}$ with values in
$\R^m$ normed by
\[
\Vert Z\Vert_{\mathcal{M}^p}=\E\biggl[\biggl(\int_0^T \vert
Z_s\vert^2\,ds\biggr)^{p/2}\biggr]^{1/p}.
\]
\end{longlist}
In the following we keep the same notation $C$ for all finite,
nonnegative constants that appear in our computations; they may depend
on known parameters deriving from assumptions and on $T$ but not on any
of the approximation and discretization parameters. In the same spirit,
we keep the same notation $\eta$ for all finite, positive constants
that we can take as small as we want independently of the approximation
and discretization parameters.

%s2.2 ###
\subsection{Some results on BMO martingales}\label{sec2.2}
In our work, the space of BMO martingales play a key role for the a
priori estimates needed in our analysis of BSDEs. We refer the reader
to \cite{Kazamaki-94} for the theory of BMO martingales and we just
recall the properties that we will use in the sequel. Let $\Phi_t=\int
_0^t \phi_s \,dW_s$, $t \in[0,T]$, be a real square integrable martingale
with respect to the Brownian filtration. Then $\Phi$ is a BMO
martingale if
\[
\Vert\Phi\Vert_{\mathrm{BMO}}=\sup_{\tau\in[0,T]} \E
[\langle\Phi
\rangle
_T-\langle\Phi\rangle_{\tau}|\mathcal{F}_{\tau}
]^{1/2}=\sup
_{\tau\in[0,T]} \E\biggl[\int_{\tau}^T \phi_s^2\,ds\big|\mathcal
{F}_{\tau
}\biggr]^{1/2} < +\infty,
\]
where the supremum is taken over all stopping times in $[0,T]$;
$\langle\Phi\rangle$ denotes the quadratic variation of $\Phi$. In
our case, the very important feature of BMO martingales is the
following lemma.
\begin{lem}
\label{proprietes martingales BMOs}
Let $\Phi$ be a BMO martingale. Then we have:
\begin{longlist}[(1)]
\item[(1)] The stochastic exponential
\[
\mathcal{E}(\Phi)_t=\mathcal{E}_t=\exp\biggl(\int_0^t \phi_s
\,dW_s-\frac
{1}{2}\int_0^t \vert\phi_s\vert^2\,ds \biggr),\qquad 0
\leq t \leq T,
\]
is a uniformly integrable martingale.
\item[(2)] Thanks to the reverse H\"{o}lder inequality, there exists
$p>1$ such that $\mathcal{E}_T \in L^p$. The maximal $p$ with this
property can be expressed in terms of the BMO norm of $\Phi$.
\item[(3)] $\forall n \in\N^*$, $\E[(\int_{0}^T
\vert\phi_s\vert^2\,ds)^n] \leq n! \Vert
\Phi\Vert_{\mathrm{BMO}}^{2n}$.
\end{longlist}
\end{lem}

%s2.3 ###
\subsection{The backward--forward system}\label{sec2.3}
Given functions $b$, $\sigma$, $g$ and $f$, for $x \in\R^d$ we will
deal with the solution $(X,Y,Z)$ to the following system of (decoupled)
backward--forward stochastic differential equations: for $t \in[0,T]$,
%
%e5 ###
%e4 ###
%
\begin{eqnarray}
\label{EDS}
X_t &=& x +\int_0^t b(s,X_s)\,ds+\int_0^t \sigma(s) \,dW_s,\\
\label{EDSR}
Y_t &=& g(X_T) +\int_t^T f(s,X_s,Y_s,Z_s)\,ds-\int_t^T Z_s \,dW_s.
\end{eqnarray}
For the functions that appear in the above system of equations we give
some general assumptions.

\subsubsection*{\textup{(HX0)}} $b \dvtx [0,T] \times\R^d \rightarrow\R^d$,
$\sigma
\dvtx [0,T] \rightarrow\R^{d\times d}$ are measurable functions. There
exist four positive constants $M_b$, $K_b$, $M_{\sigma}$ and
$K_{\sigma
}$ such that $\forall t,t' \in[0,T]$, $\forall x,x' \in\R^d$,
\begin{eqnarray*}
\vert b(t,x)\vert& \leq& M_b(1+\vert x
\vert),\\
\vert b(t,x)-b(t',x')\vert&\leq& K_b (\vert
x-x'\vert+\vert t-t'\vert^{1/2}),\\
\vert\sigma(t)\vert&\leq& M_{\sigma},\\
\vert\sigma(t)-\sigma(t')\vert&\leq& K_{\sigma}
\vert t-t'\vert.
\end{eqnarray*}

\subsubsection*{\textup{(HY0)}} $f \dvtx [0,T] \times\R^d \times\R\times\R^{1
\times d} \rightarrow\R$, $g \dvtx \R^d \rightarrow\R$ are measurable
functions. There exist five positive constants $M_f$, $K_{f,x}$,
$K_{f,y}$, $K_{f,z}$ and $M_g$ such that $\forall t \in[0,T]$,
$\forall x,x' \in\R^d$, $\forall y,y' \in\R$, $\forall z,z' \in\R
^{1\times d}$,
\begin{eqnarray*}
\vert f(t,x,y,z)\vert&\leq& M_f(1+\vert y
\vert+\vert z\vert^2),\\
\vert f(t,x,y,z)-f(t,x',y',z')\vert&\leq& K_{f,x}
\vert x-x'\vert+K_{f,y}\vert y-y'\vert\\
&&{}+\bigl(K_{f,z}+L_{f,z}(\vert z\vert+\vert z'
\vert)\bigr)\vert z-z'\vert,\\
\vert g(x)\vert&\leq& M_g.
\end{eqnarray*}
We next recall some results on BSDEs with quadratic growth. For their
original version and their proof we refer to
\cite{Kobylanski-00,Briand-Confortola-08} and \cite{Imkeller-dosReis-09}.
\begin{theorem}
\label{norme BMO, def r q}
Under \textup{(HX0)}, \textup{(HY0)}, the system (\ref{EDS})--(\ref{EDSR}) has
a unique solution $(X,Y,Z) \in\mathcal{S}^2 \times\mathcal
{S}^{\infty
} \times\mathcal{M}^2$. The martingale $Z\ast W$ belongs to the space
of BMO martingales and $\Vert Z \ast W\Vert_{\mathrm{BMO}}$ only
depends on $T$,
$M_g$ and $M_f$. Moreover, there exists $r>1$ such that $\mathcal{E}(Z
\ast W) \in L^r$.
\end{theorem}

%s3 ###
\section{Some useful estimates of $Z$}\label{sec3}

%s3.1 ###
\subsection{A first bound for $Z$}\label{sec3.1}
\begin{theorem}
\label{borne globale sur Z}
Suppose that \textup{(HX0)}, \textup{(HY0)} hold and that $g$ is Lipschitz with
Lipschitz constant $K_{g}$. Then, there exists a version of $Z$ such
that, $\forall t \in[0,T]$,
\[
\vert Z_t\vert\leq e^{(2K_b+K_{f,y})T}M_{\sigma}(K_g+TK_{f,x}).
\]
\end{theorem}
\begin{pf}
First, we suppose that $b$, $g$ and $f$ are differentiable with respect
to $x$, $y$ and $z$. Then $(X,Y,Z)$ is differentiable with respect to
$x$ and $(\nabla X,\nabla Y,\nabla Z)$ is the solution of
%
%e7 ###
%e6 ###
%
\begin{eqnarray}
\label{EDS nabla X}
\nabla X_t &=& I_d+\int_0^t \nabla b(s,X_s)\nabla X_s \,ds,\\
\label{EDSR nabla Y nabla Z}
\nabla Y_t &=& \nabla g(X_T)\nabla X_T -\int_t^T \nabla Z_s \,dW_s \nonumber\\
&&{}+\int_t^T \nabla_x f(s,X_s,Y_s,Z_s)\nabla X_s+\nabla_y
f(s,X_s,Y_s,Z_s)\nabla Y_s\,ds \\
&&{}+\int_t^T\nabla_z f(s,X_s,Y_s,Z_s)\nabla Z_s\,ds, \nonumber
\end{eqnarray}
where $\nabla X_t=(\partial X^i_t/\partial x^j)_{1 \leq i,j
\leq d}$, $\nabla Y_t=\supmathitt{(\partial Y_t/\partial x^j)_{1
\leq j \leq d}} \in\R^{1 \times d}$, $\nabla Z_t=(\partial
Z^i_t/\break
\partial x^j)_{1 \leq i,j \leq d}$ and $\int_t^T \nabla
Z_s \,dW_s$ means
\[
\sum_{1 \leq i \leq d} \int_t^T (\nabla Z_s)^i \,dW^i_s
\]
with $(\nabla Z)^i$ denoting the $i$th line of the $d\times d$ matrix
process $\nabla Z$. Thanks to usual transformations on the BSDE we obtain
\begin{eqnarray*}
&&e^{\int_0^t\nabla_y f(s,X_s,Y_s,Z_s)\,ds}\nabla Y_t\\
&&\qquad= e^{\int_0^T\nabla_y f(s,X_s,Y_s,Z_s)\,ds}\nabla g(X_T)\nabla X_T\\
&&\qquad\quad{}-\int_t^T e^{\int_0^s\nabla_y f(u,X_u,Y_u,Z_u)\,du} \nabla Z_s
\,d\tilde
{W}_s\\
&&\qquad\quad{} + \int_t^T e^{\int_0^s\nabla_y f(u,X_u,Y_u,Z_u)\,du}\nabla_x
f(s,X_s,Y_s,Z_s)\nabla X_s \,ds
\end{eqnarray*}
with $d\tilde{W}_s = dW_s-\nabla_z f(s,X_s,Y_s,Z_s)\,ds$. We have
\begin{eqnarray*}
&&\biggl\Vert\int_0^{\cdot} \nabla_zf(s,X_s,Y_s,Z_s)\,dW_s\biggr\Vert
^2_{\mathrm{BMO}}\\
&&\qquad = \sup_{\tau\in[0,T]} \E\biggl[\int_{\tau}^T\vert
\nabla_z f(s,X_s,Y_s,Z_s)\vert^2\,ds \big| \mathcal{F}_{\tau}
\biggr]\\
&&\qquad\leq C\biggl(1+\sup_{\tau\in[0,T]} \E\biggl[\int_{\tau}^T
\vert Z_s\vert^2 \,ds\big| \mathcal{F}_{\tau} \biggr]
\biggr)\\
&&\qquad=C (1+\Vert Z\ast W\Vert^2_{\mathrm{BMO}}).
\end{eqnarray*}
Since $Z\ast W$ belongs to the space of BMO martingales,
\[
\biggl\Vert\int_0^{\cdot} \nabla_zf(s,X_s,Y_s,Z_s)\,dW_s\biggr\Vert
_{\mathrm{BMO}}<+\infty.
\]
Lemma \ref{proprietes martingales BMOs} gives us that $\mathcal
{E}(\int
_0^{\cdot} \nabla_zf(s,X_s,Y_s,Z_s)\,dW_s)_t$ is a uniformly integrable
martingale so we are able to apply Girsanov's theorem: there exists a
probability $\mathbb{Q}$ under which $(\tilde{W})_{t \in[0,T]}$ is a
Brownian motion. Then,
\begin{eqnarray*}
&&e^{\int_0^t\nabla_y f(s,X_s,Y_s,Z_s)\,ds}\nabla Y_t\\
&&\qquad= \E^{\Q}\biggl[e^{\int_0^T\nabla_y f(s,X_s,Y_s,Z_s)\,ds}\nabla
g(X_T)\nabla X_T\\
&&\hspace*{16.8pt}\qquad\quad{}+ \int_t^T e^{\int_0^s\nabla_y f(u,X_u,Y_u,Z_u)\,du}\nabla_x
f(s,X_s,Y_s,Z_s)\nabla X_s\,ds \big| \mathcal{F}_t \biggr]
\end{eqnarray*}
and
%
%e8 ###
%
\begin{equation}
\label{majoration nabla Y}
\vert\nabla Y_t\vert\leq e^{(K_b+K_{f,y})T}(K_g+TK_{f,x}),
\end{equation}
because $\vert\nabla X_t\vert\leq e^{K_b T}$. Moreover,
thanks to the
Malliavin calculus, it is classical to show that a version of $(Z_t)_{t
\in[0,T]}$ is given by $(\nabla Y_t(\nabla X_t)^{-1} \sigma(t))_{t
\in
[0,T]}$. So we obtain
\[
\vert Z_t\vert\leq e^{K_bT}M_{\sigma}\vert\nabla
Y_t\vert\leq
e^{(2K_b+K_{f,y})T}M_{\sigma}(K_g+TK_{f,x}) \qquad\mbox{a.s.},
\]
because $\vert\nabla X_t^{-1}\vert\leq e^{K_b T}$.

When $b$, $g$ and $f$ are not differentiable, we can also prove the
result by a standard approximation and stability results for BSDEs with
linear growth.
\end{pf}
\begin{rem}
Thanks to Theorem \ref{borne globale sur Z}, the generator $f$
becomes a Lipschitz function with respect to $z$, so we are able to use
standard results about time discretization of BSDEs. In this case, we
obtain that the error of approximation is lower than $Cn^{-1}$ with $n$
the number of discretization times (see, e.g., \cite
{Bouchard-Touzi-04,Gobet-Lemor-Warin-05}). Let us notice that, in all
studies about discretization of BSDEs, we do not care about the
constant in the error bound; we only consider the asymptotic speed of
convergence. But, with a practical point of view, the constant could
play an important role, particularly for small $n$. In our case, the
generator $f$ may be viewed as Lipschitz in $z$ with a Lipschitz
constant $Ce^{(2K_b+K_{f,y})T}$. So, if we apply the standard result,
the generic constant in the rate of convergence will be in the order of
$Ce^{Ce^{2(2K_b+K_{f,y})T}}$. This is, of course, not desirable because
it blows up when $K_b$, $K_{f,y}$ or $T$ increase. We think that it
could be interesting to see if we are able to observe such a phenomena
with numerical simulation.
\end{rem}

%s3.2 ###
\subsection{A time dependent estimate of $Z$}\label{sec3.2}
We will introduce two alternative assumptions.

\subsubsection*{\textup{(HX1)}}
$b$ is differentiable with respect to $x$ and $\sigma$ is
differentiable with respect to $t$. There exists $\lambda\in\R^+$
such that
$\forall\eta\in\mathbb{R}^d$
%
%e9 ###
%
\begin{equation}
\label{hypothese sur nabla b}
\vert\supmathitt{\eta}\sigma(s)[\supmathitt{\sigma
(s)}\supmathitt{\nabla b(s,x)}
-\supmathitt{\sigma'(s)}]\eta
\vert
\leq\lambda\vert\supmathitt{\eta}\sigma(s)\vert^2.
\end{equation}

\subsubsection*{\textup{(HX1$'$)}}
$\sigma$ is invertible and $\forall t \in[0,T]$, $\vert\sigma
(t)^{-1}\vert\leq M_{\sigma^{-1}}$.

\subsubsection*{Example}
Assumption (HX1) is verified when, $\forall s \in[0,T]$, $\nabla
b(s,\cdot)$ commutes with $\sigma(s)$ and $\exists A \dvtx [0,T] \rightarrow
\R
^{d \times d}$ bounded such that $\sigma'(t)=\sigma(t)A(t)$.
\begin{theorem}
\label{borne temporelle sur Z}
Suppose that \textup{(HX0)}, \textup{(HY0)} hold and that \textup{(HX1)} or
\textup{(HX1$'$)} holds. Moreover, suppose that $g$ is lower (or upper)
semi-continuous. Then there exists a version of $Z$ and there exist two
constants $C,C' \in\R^+$ that depend only in $T$, $M_g$, $M_f$,
$K_{f,x}$, $K_{f,y}$, $K_{f,z}$ and $L_{f,z}$ such that, $\forall t \in[0,T[$,
\[
\vert Z_t\vert\leq C+C'(T-t)^{-1/2}.
\]
\end{theorem}
\begin{pf}
In a first time, we will suppose that (HX1) holds and that $f$, $g$ are
differentiable with respect to $x$, $y$ and $z$. Then $(Y,Z)$ is
differentiable with respect to $x$ and $(\nabla Y,\nabla Z)$ is the
solution of the BSDE
\begin{eqnarray*}
\nabla Y_t &=& \nabla g(X_T)\nabla X_T - \int_t^T \nabla Z_s \,dW_s\\
&&{} +\int_t^T \nabla_x f(s,X_s,Y_s,Z_s) \nabla X_s + \nabla_y
f(s,X_s,Y_s,Z_s) \nabla Y_s\,ds\\
&&{}+ \int_t^T\nabla_z f(s,X_s,Y_s,Z_s) \nabla Z_s \,ds.
\end{eqnarray*}
Thanks to usual transformations we obtain
\begin{eqnarray*}
& & e^{\int_0^t\nabla_y f(s,X_s,Y_s,Z_s)\,ds}\nabla Y_t\\
&&\quad{}+\int_0^t e^{\int_0^s\nabla_y f(u,X_u,Y_u,Z_u)\,du}\nabla_x
f(s,X_s,Y_s,Z_s)\nabla X_s \,ds\\
& &\qquad= e^{\int_0^T\nabla_y f(s,X_s,Y_s,Z_s)\,ds}\nabla g(X_T)\nabla X_T\\
&&\qquad\quad{}+ \int_0^T e^{\int_0^s\nabla_y f(u,X_u,Y_u,Z_u)\,du}\nabla_x
f(s,X_s,Y_s,Z_s)\nabla X_s \,ds\\
&&\qquad\quad{} -\int_t^T e^{\int_0^s\nabla_y f(u,X_u,Y_u,Z_u)\,du} \nabla Z_s
\,d\tilde{W}_s
\end{eqnarray*}
with $d\tilde{W}_s = dW_s-\nabla_z f(s,X_s,Y_s,Z_s)\,ds$. We can rewrite
it as
%
%e10 ###
%
\begin{equation}
F_t= F_T-\int_t^T e^{\int_0^s\nabla_y f(u,X_u,Y_u,Z_u)\,du} \nabla Z_s
\,d\tilde{W}_s
\end{equation}
with
\begin{eqnarray*}
F_t &:=& e^{\int_0^t\nabla_y f(s,X_s,Y_s,Z_s)\,ds}\nabla Y_t\\
&&{}+\int_0^t e^{\int_0^s\nabla_y f(u,X_u,Y_u,Z_u)\,du}\nabla_x
f(s,X_s,Y_s,Z_s)\nabla X_s \,ds.
\end{eqnarray*}
$Z\ast W$ belongs to the space of BMO martingales so we are able to
apply Girsanov's theorem: there exists a probability $\mathbb{Q}$ under
which $(\tilde{W})_{t \in[0,T]}$ is a Brownian motion. Thanks to the
Malliavin calculus, it is possible to show that $(\nabla Y_t(\nabla
X_t)^{-1} \sigma(t))_{t \in[0,T]}$ is a version of $Z$. Now we define
\begin{eqnarray*}
\alpha_t &:= & \int_0^t e^{\int_0^s\nabla_y
f(u,X_u,Y_u,Z_u)\,du}\nabla_x
f(s,X_s,Y_s,Z_s)\nabla X_s \,ds\, (\nabla X_t)^{-1}\sigma(t),\\
\tilde{Z}_t &:=& F_t (\nabla X_t)^{-1} \sigma(t)= e^{\int_0^t\nabla_y
f(s,X_s,Y_s,Z_s)\,ds} Z_t +\alpha_t \qquad\mbox{a.s.},\\
\tilde{F}_t &:=& e^{\lambda t} F_t (\nabla X_t)^{-1}.
\end{eqnarray*}
Since $d\nabla X_t=\nabla b(t,X_t)\nabla X_t \,dt$, then $d(\nabla
X_t)^{-1} = - (\nabla X_t)^{-1}\nabla b(t,X_t)\,dt$ and thanks to It\^
{o}'s formula,
\[
d\tilde{Z}_t=dF_t(\nabla X_t)^{-1} \sigma(t)-F_t (\nabla
X_t)^{-1}\nabla b(t,X_t)\sigma(t)\,dt+F_t(\nabla X_t)^{-1}\sigma'(t)\,dt
\]
and
\[
d(e^{\lambda t}\tilde{Z}_t)=\tilde{F}_t\bigl(\lambda Id-\nabla
b(t,X_t)\bigr)\sigma(t)\,dt+\tilde{F}_t\sigma'(t)\,dt+e^{\lambda t}
\,dF_t(\nabla
X_t)^{-1} \sigma(t).
\]
Finally,
\begin{eqnarray*}
d\vert e^{\lambda t} \tilde{Z}_t\vert^2&=&2
\bigl[\lambda\vert\tilde{F}_t\sigma(t)\vert^2-\tilde
{F}_t\sigma(t)[\supmathitt{\sigma
(t)}\supmathitt{\nabla b(t,X_t)}-\supmathitt{\sigma
'(t)}]\supmathitt{\tilde{F}_t}\bigr]\,dt\\
&&{}+d\langle M \rangle_t+dM_t^*
\end{eqnarray*}
with $M_t:=\int_0^t e^{\lambda s} \,dF_s(\nabla X_s)^{-1}\sigma(s)$ and
$M_t^*$ a $\mathbb{Q}$-martingale. Thanks to the assumption (HX1) we
are able to conclude that $\vert e^{\lambda t} \tilde{Z}_t
\vert^2$ is a
$\mathbb{Q}$-submartingale. Hence,
\begin{eqnarray*}
&& \mathbb{E}^{\mathbb{Q}} \biggl[\int_t^T e^{2\lambda s}\vert
\tilde{Z}_s\vert^2\,ds \big|\mathcal{F}_t\biggr]\\
&&\qquad \geq e^{2\lambda t}\vert\tilde{Z}_t\vert^2(T-t)\\
&&\qquad\geq e^{2\lambda t}\bigl\vert e^{\int_0^t\nabla_y
f(s,X_s,Y_s,Z_s)\,ds}Z_t+\alpha_t\bigr\vert^2(T-t) \qquad\mbox{a.s.},
\end{eqnarray*}
which implies
\begin{eqnarray*}
\vert Z_t\vert^2(T-t) & = & e^{-2\lambda t}e^{-2\int
_0^t\nabla_y
f(s,X_s,Y_s,Z_s)\,ds}e^{2\lambda t}\\
&&{} \times\bigl\vert e^{\int_0^t\nabla_y
f(s,X_s,Y_s,Z_s)\,ds}Z_t+\alpha_t-\alpha_t\bigr\vert^2(T-t)\\
&\leq& C\bigl(e^{2\lambda t}\bigl\vert e^{\int_0^t\nabla_y
f(s,X_s,Y_s,Z_s)\,ds}Z_t+\alpha_t\bigr\vert^2+1\bigr)(T-t)\\
&\leq& C\biggl( \mathbb{E}^{\mathbb{Q}} \biggl[\int_t^T e^{2\lambda
s}\vert\tilde{Z}_s\vert^2\,ds \big|\mathcal{F}_t
\biggr] +(T-t)\biggr)\qquad
\mbox{a.s.}
\end{eqnarray*}
with $C$ a constant that only depends on $T$, $K_b$, $M_{\sigma}$,
$K_{f,x}$, $K_{f,y}$ and $\lambda$. Moreover, we have, a.s.,
\begin{eqnarray*}
\mathbb{E}^{\mathbb{Q}} \biggl[\int_t^T e^{2\lambda s}\vert
\tilde{Z}_s\vert^2\,ds \big|\mathcal{F}_t\biggr] & \leq& C
\mathbb{E}^{\mathbb
{Q}} \biggl[\int_t^T \vert Z_s\vert^2+\vert\alpha
_s\vert^2\,ds \big|\mathcal
{F}_t\biggr]\\
& \leq& C\bigl(\Vert Z\Vert^2_{\mathrm{BMO}(\mathbb{Q})}
+(T-t)\bigr).
\end{eqnarray*}
But $\Vert Z\Vert_{\mathrm{BMO}(\mathbb{Q})}$ does not depend on
$K_g$ because
$(Y,Z)$ is a solution of the following quadratic BSDE:
%
%e11 ###
%
\begin{eqnarray}
\label{EDSR modifiee sous Q}
Y_t&=&g(X_T)+\int_t^T\bigl(f(s,X_s,Y_s,Z_s)-Z_s\nabla
_zf(s,X_s,Y_s,Z_s)\bigr)\,ds\nonumber\\[-8pt]\\[-8pt]
&&{}-\int_t^TZ_s\,d\tilde{W}_s.\nonumber
\end{eqnarray}
Finally, $\vert Z_t\vert\leq C
(1+(T-t)^{-1/2})$ a.s.

When $\sigma$ is invertible, the inequality (\ref{hypothese sur nabla
b}) is verified with $\lambda:=M_{\sigma^{-1}}(M_{\sigma
}K_b+K_{\sigma
})$. Since this $\lambda$ does not depend on $\nabla b$ and $\sigma'$,
we can prove the result when $b(t,\cdot)$ and $\sigma$ are not
differentiable by a standard approximation and stability results for
BSDEs with linear growth. So, we are allowed to replace assumption
(HX1) by (HX1$'$).

When $f$ is not differentiable and $g$ is only Lipschitz, we can prove
the result by a standard approximation and stability results for linear
BSDEs. But we notice that our estimation on $Z$ does not depend on
$K_g$. This allows us to weaken the hypothesis on $g$ further; when $g$
is only lower or upper semi-continuous the result stays true. The proof
is the same as the proof of Proposition 4.3 in \cite{Delbaen-Hu-Bao-09}.
\end{pf}
\begin{rem}
\label{estimation temporelle plus precise pour Z}
The previous proof gives us a more precise estimation for a version of
$Z$ when $f$ is differentiable with respect to $z\dvtx\forall t \in[0,T]$,
\[
\vert Z_t\vert\leq C+C'\E^{\Q}\biggl[\int_t^T
\vert Z_s\vert^2\,ds
\big|\mathcal{F}_t\biggr]^{1/2}(T-t)^{-1/2}.
\]
\end{rem}
\begin{rem}
When assumption (HX1) or (HX1$'$) is not verified, the process $Z$
may blow up before $T$. Zhang gives an example of such a phenomenon in
dimension 1; we refer the reader to Example 1 in \cite{Zhang-05}.
\end{rem}

%s3.3 ###
\subsection{Zhang's path regularity theorem}\label{sec3.3}
Let $0=t_0<t_1<\cdots< t_n=T$ be any given partition of $[0,T]$ and
denote $\delta_n$ the mesh size of this partition. We define a set of
random variables
\[
\bar{Z}_{t_i} = \frac{1}{t_{i+1}-t_i}\E\biggl[\int_{t_i}^{t_{i+1}} Z_s\,ds
\big|\mathcal{F}_{t_i}\biggr]\qquad \forall i \in\{0,\ldots
,n-1\}.
\]
Then we are able to give a more detailed version at
Theorem 3.4.3 in \cite{Zhang-01}.
\begin{theorem}
\label{thm de Zhang}
Suppose that \textup{(HX0)}, \textup{(HY0)} hold and $g$ is a Lipschitz function
with Lipschitz constant $K_{g}$. Then we have
\[
\sum_{i=0}^{n-1} \E\biggl[ \int_{t_i}^{t_{i+1}} \vert
Z_t-\bar{Z}_{t_i}\vert^2\,dt\biggr] \leq C(1+K_{g}^2)\delta_n,
\]
where C is a positive constant independent of $\delta_n$ and $K_{g}$.
\end{theorem}
\begin{pf}
We will follow the proof of Theorem 5.6 in \cite{Imkeller-dosReis-09};
we just need to specify how the estimate depends on $K_g$. First, it is
not difficult to show that $\bar{Z}_{t_i}$ is the best $\mathcal
{F}_{t_i}$-measurable approximation of $Z$ in $\mathcal
{M}^2([t_i,t_{i+1}])$, that is,
\[
\E\biggl[\int_{t_i}^{t_{i+1}} \vert Z_t-\bar{Z}_{t_i}
\vert^2\,dt\biggr] =
\inf
_{Z_i \in L^2(\Omega,\mathcal{F}_{t_i})} \E\biggl[\int_{t_i}^{t_{i+1}}
\vert Z_t-Z_{i}\vert^2\,dt\biggr].
\]
In particular,
\[
\E\biggl[\int_{t_i}^{t_{i+1}} \vert Z_t-\bar{Z}_{t_i}
\vert^2\,dt\biggr]
\leq\E\biggl[\int_{t_i}^{t_{i+1}} \vert Z_t-Z_{t_i}
\vert^2\,dt\biggr].
\]
In the same spirit as previous proofs, we suppose in a first time that
$b$, $g$ and $f$ are differentiable with respect to $x$, $y$ and $z$. So,
\[
Z_t-Z_{t_i}=\nabla Y_t (\nabla X_t)^{-1} \sigma(t)-\nabla Y_{t_i}
(\nabla X_{t_i})^{-1} \sigma(t_i)=I_1+I_2+I_3 \qquad\mbox{a.s.},
\]
with $I_1= \nabla Y_t (\nabla X_t)^{-1} (\sigma(t)-\sigma(t_i))$,
$I_2=\nabla Y_t ((\nabla X_t)^{-1}-(\nabla X_{t_i})^{-1}) \sigma(t_i)$
and $I_3=\nabla(Y_t-Y_{t_i}) (\nabla X_{t_i})^{-1} \sigma(t_i)$.
First, thanks to the estimation (\ref{majoration nabla Y}) we have
\[
\vert I_1\vert^2 \leq\vert\nabla Y_t\vert
^2e^{2K_bT} K_{\sigma}^2 \vert t_{i+1}-t_i\vert^2 \leq
C(1+K_g^2)\delta_n^2.
\]
We obtain the same estimation for $\vert I_2\vert$ because
\[
\vert(\nabla X_t)^{-1}-(\nabla X_{t_i})^{-1}\vert\leq
\biggl\vert\int_{t_i}^t (\nabla X_s)^{-1}\nabla b(s,X_s)\,ds
\biggr\vert\leq K_b e^{K_b T}\vert t-t_i\vert.
\]
Last, $\vert I_3\vert\leq M_{\sigma} e^{K_bT}
\vert\nabla Y_t-\nabla Y_{t_i}\vert$. So,
\[
\sum_{i=0}^{n-1} \E\biggl[\int_{t_i}^{t_{i+1}} \vert I_3
\vert^2 \,dt\biggr]
\leq C\delta_n \sum_{i=0}^{n-1} \E\Bigl[ \esssup_{t \in
[t_i,t_{i+1}]} \vert\nabla Y_t-\nabla Y_{t_i}\vert
^2\Bigr].
\]
By using the BSDE (\ref{EDSR nabla Y nabla Z}), (HY0), the estimate on
$\nabla X_s$ and the estimate (\ref{majoration nabla Y}), we have
\begin{eqnarray*}
\vert\nabla Y_t-\nabla Y_{t_i}\vert^2
&\leq& C\biggl( \int_{t_i}^{t} \bigl(C(1+K_g)+\vert\nabla_z
f(s,X_s,Y_s,Z_s)\vert\vert\nabla Z_s\vert
\bigr)\,ds\biggr)^2\\
&&{}+C\biggl( \int_{t_i}^t \nabla Z_s\,dW_s \biggr)^2.
\end{eqnarray*}
The inequalities of H\"{o}lder and Burkholder--Davis--Gundy give us
\begin{eqnarray*}
&&\sum_{i=0}^{n-1} \E\Bigl[ \esssup_{t \in[t_i,t_{i+1}]}
\vert\nabla Y_t-\nabla Y_{t_i}\vert^2\Bigr]\\
&&\qquad\leq C(1+K_g^2)+C\sum_{i=0}^{n-1} \E\biggl( \int
_{t_i}^{t_{i+1}}\vert\nabla_z f(s,X_s,Y_s,Z_s)\vert
\vert\nabla Z_s\vert \,ds
\biggr)^2\\
&&\qquad\quad{}+C\E\biggl( \int_{t_i}^{t_{i+1}} \vert\nabla Z_s\vert
^2\,ds \biggr)\\
&&\qquad\leq C(1+K_g^2)\\
&&\qquad\quad{}+C\E\biggl[ \biggl( \int_{0}^{T}\vert\nabla
_z f(s,X_s,Y_s,Z_s)\vert\vert\nabla Z_s\vert
\,ds\biggr)^2+ \int_{0}^{T} \vert\nabla Z_s\vert
^2\,ds\biggr]\\
&&\qquad\leq C(1+K_g^2)\\
&&\qquad\quad{}+C\E\biggl[ \biggl(\int_{0}^{T}(1+\vert
Z_s\vert^2)\,ds\biggr)\biggl(\int_0^T\vert\nabla
Z_s\vert^2\,ds\biggr)
+ \int
_{0}^{T} \vert\nabla Z_s\vert^2\,ds\biggr]\\
&&\qquad\leq C(1+K_g^2)\\
&&\qquad\quad{}+C\biggl(1+\E\biggl[ \biggl(\int_{0}^{T}\vert
Z_s\vert^2\,ds\biggr)^{p}\biggr]^{1/p}\biggr)\E\biggl[
\biggl(\int
_0^T\vert\nabla Z_s\vert^2\,ds\biggr)^{q}\biggr]^{1/q}
\end{eqnarray*}
for all $p>1$ and $q>1$ such that $1/p+1/q=1$. But, $(\nabla Y, \nabla
Z)$ is the solution of BSDE (\ref{EDSR nabla Y nabla Z}) so, from
Corollary 9 in \cite{Briand-Confortola-08}, there exists $q$ that only
depends on $\Vert Z\ast W\Vert_{\mathrm{BMO}}$ such that
\[
\E\biggl[\biggl(\int_{0}^{T} \vert\nabla Z_s\vert
^2\,ds\biggr)^{q}\biggr]^{1/q}
\leq C(1+K_g^2).
\]
Moreover, we can apply Lemma \ref{proprietes martingales BMOs} to
obtain the estimate
\[
\E\biggl[ \biggl(\int_{0}^{T}\vert Z_s\vert^2\,ds
\biggr)^{p}\biggr]^{1/p}
\leq C\Vert Z\Vert_{\mathrm{BMO}}^{2} \leq C.
\]
Finally,
\[
\sum_{i=0}^{n-1} \E\biggl[\int_{t_i}^{t_{i+1}} \vert I_3
\vert^2 \,dt\biggr]
\leq C(1+K_g^2)\delta_n
\]
and
\begin{eqnarray*}
\sum_{i=0}^{n-1}\E\biggl[\int_{t_i}^{t_{i+1}} \vert Z_t-\bar
{Z}_{t_i}\vert^2\,dt\biggr] &\leq& \sum_{i=0}^{n-1} \E
\biggl[\int
_{t_i}^{t_{i+1}} (\vert I_1\vert^2+\vert
I_2\vert^2+\vert I_3\vert^2
)\,dt\biggr]\\
& \leq& C(1+K_g^2)\delta_n.
\end{eqnarray*}
\upqed\end{pf}

%s4 ###
\section{Convergence of a modified time discretization scheme for the BSDE}\label{sec4}
%s4.1 ###
\subsection{An approximation of the quadratic BSDE}\label{sec4.1}
In a first time we will approximate our quadratic BSDE (\ref{EDSR}) by
another one. We set \mbox{$\eps\in\ ]0,T[$} and $N \in\N$. Let $(Y_t^{N,\eps
},Z_t^{N,\eps})$ be the solution of the BSDE
%
%e12 ###
%
\begin{equation}
\label{EDSR_modifiee}
Y_t^{N,\eps} = g_N(X_T)+\int_t^T f^{\eps}(s,X_s,Y_s^{N,\eps
},Z_s^{N,\eps
})\,ds-\int_t^T Z_s^{N,\eps} \,dW_s
\end{equation}
with
\[
f^{\eps}(s,x,y,z) := \mathbh{1}_{s\leq T-\eps} f(s,x,y,z)+ \mathbh
{1}_{s > T-\eps} f(s,x,y,0)
\]
and $g_N$ a Lipschitz approximation of $g$ with Lipschitz constant $N$.
$f^{\eps}$ verifies assumption (HY0) with the same constants as $f$.
Since $g_N$ is a Lipschitz function, $Z^{N,\eps}$ has a bounded version
and the BSDE (\ref{EDSR_modifiee}) is a BSDE with a linear growth.
Moreover, we can apply Theorem \ref{borne temporelle sur Z} to obtain
the following proposition.
\begin{prop}
\label{majoration Z modifie} Let us assume that \textup{(HX0)},
\textup{(HY0)} and \textup{(HX1)} or \textup{(HX1$'$)} hold. There
exists a version of $Z^{N,\eps}$ and there exist three constants
$M_{z,1},M_{z,2}$, $M_{z,3} \in\R^+$ that do not depend on $N$ and $\eps$
such that, $\forall s \in[0,T]$,
\[
\vert Z_s^{N,\eps}\vert\leq\biggl( M_{z,1}+\frac
{M_{z,2}}{(T-s)^{1/2}}\biggr) \wedge\bigl(M_{z,3}(N+1)\bigr).
\]
\end{prop}

Thanks to BMO tools we have a stability result for quadratic BSDEs (see
\cite{Briand-Confortola-08} and \cite{Imkeller-dosReis-09}).
\begin{prop}
\label{erreur premiere approximation EDSR}
Let us assume that \textup{(HX0)} and \textup{(HY0)} hold. There exists a
constant $C$ that does not depend on $N$ and $\eps$ such that
\[
\E\Bigl[\sup_{t \in[0,T]}\vert Y_t^{N,\eps}-Y_t\vert
^2\Bigr]+ \E
\biggl[
\int_0^T \vert Z_t^{N,\eps}-Z_t\vert^2\,dt \biggr] \leq C
\bigl(e_1(N)+e_2(N,\eps)\bigr)
\]
with
\begin{eqnarray*}
e_1(N)&:=& \E[\vert g_N(X_T)-g(X_T)\vert^{2q}]^{1/q},
\\
e_2(N,\eps)&:=& \E\biggl[ \biggl(\int_{T-\eps}^T \vert
f(t,X_t,Y_t^{N,\eps},Z_t^{N,\eps})-f(t,X_t,Y_t^{N,\eps},0)
\vert \,dt \biggr)^{2q}\biggr]^{1/q}
\end{eqnarray*}
and $q$ defined in Theorem \ref{norme BMO, def r q}.
\end{prop}
\begin{rem}
The authors of \cite{Imkeller-dosReis-09} obtain this result with
$q^2$ instead of $q$. Nevertheless, we are able to obtain the good
result by applying the estimates of~\cite{Briand-Confortola-08}.
\end{rem}

Then, in a second time, we will approximate our modified
backward--forward system by a discrete-time one. We will slightly modify
the classical discretization by using a nonequidistant net with $2n+1$
discretization times. We define the $n+1$ first discretization times on
$[0,T-\eps]$ by
\[
t_{k}=T\biggl( 1-\biggl( \frac{\eps}{T}\biggr)^{k/n} \biggr)
\]
and we use an equidistant net on $[T-\eps,T]$ for the last $n$
discretization times
\[
t_k=T-\biggl( \frac{2n-k}{n} \biggr)\eps,\qquad n \leq k \leq2n.
\]
We denote the time step by $(h_k:=t_{k+1}-t_k)_{0\leq k \leq2n-1}$.
We consider $(X^n_{t_k})_{0\leq k \leq2n}$ the classical
Euler scheme for $X$ given by
%
%e13 ###
%
\begin{eqnarray}
X^n_0 &=& x, \nonumber\\[-8pt]\\[-8pt]
X^n_{t_{k+1}}&=&X^n_{t_k}+h_kb(t_k,X^n_{t_k})+\sigma(t_k)(W_{t_{k+1}}-W_{t_k})\nonumber
\end{eqnarray}
for $0 \leq k \leq2n-1$. We denote $\rho_s \dvtx \R^{1 \times d}
\rightarrow\R^{1 \times d}$ the projection on the ball
\[
B\biggl(0,M_{z,1}+\frac{M_{z,2}}{(T-s)^{1/2}}\biggr)
\]
with $M_{z,1}$ and $M_{z,2}$ given by Proposition \ref{majoration Z modifie}.
Finally, we denote $(Y^{N,\eps,n},Z^{N,\eps,n})$ our time approximation
of $(Y^{N,\eps},Z^{N,\eps})$. This couple is obtained by a slight
modification of the classical dynamic programming equation
%
%e15 ###
%e14 ###
%
\begin{eqnarray}
Y^{N,\eps,n}_{t_{2n}} & = & g_N(X^n_{t_{2n}}), \nonumber\\
\label{Z discretisation temp}
Z^{N,\eps,n}_{t_k} &=& \rho_{t_{k+1}}\biggl(\frac{1}{h_k}\E
_{t_k}[Y^{N,\eps,n}_{t_{k+1}} (W_{t_{k+1}}-W_{t_k})]\biggr),\\
\label{Y discretisation temp}
Y^{N,\eps,n}_{t_k} & = & \E_{t_k}[Y^{N,\eps,n}_{t_{k+1}}]+h_k\E
_{t_k}[f^{\eps}(t_k,X^n_{t_k},Y^{N,\eps,n}_{t_{k+1}},Z^{N,\eps,n}_{t_k})],
\end{eqnarray}
where $0 \leq k \leq2n-1$ and $\E_{t_k}$ stand for the
conditional expectation given $\mathcal{F}_{t_k}$. Let us notice that
the classical dynamic programming equation does not use a projection in
(\ref{Z discretisation temp}); it is the only difference with our time
approximation (see, e.g., \cite{Gobet-Lemor-Warin-05} for the classical
case). This projection comes directly from the estimate of $Z$ in
Proposition \ref{majoration Z modifie}. The aim of our work is to study
the error of discretization
\[
e(N,\eps,n):= \sup_{0 \leq k \leq2n}\E[\vert
Y_{t_k}^{N,\eps,n}-Y_{t_k}\vert^2]+ \sum
_{k=0}^{2n-1}\E\biggl[
\int
_{t_k}^{t_{k+1}} \vert Z_{t_k}^{N,\eps,n}-Z_t\vert^2\,dt
\biggr].
\]
It is easy to see that
\[
e(N,\eps,n) \leq C\bigl(e_1(N)+e_2(N,\eps)+e_3(N,\eps,n)\bigr)
\]
with $e_1(N)$ and $e_2(N,\eps)$ defined in Proposition \ref{erreur
premiere approximation EDSR} and
\[
e_3(N,\eps,n):= \sup_{0 \leq k \leq2n}\E[\vert
Y_{t_k}^{N,\eps,n}-Y_{t_k}^{N,\eps}\vert^2]+ \sum
_{k=0}^{2n-1}\E
\biggl[ \int_{t_k}^{t_{k+1}} \vert Z_{t_k}^{N,\eps
,n}-Z_t^{N,\eps}\vert^2\,dt
\biggr].
\]

%s4.2 ###
\subsection{Study of the time approximation error $e_3(N,\eps,n)$}\label{sec4.2}

We need an extra assumption.

\subsubsection*{\textup{(HY1)}}
There exists a positive constant $K_{f,t}$ such that $\forall t,t' \in
[0,T]$, $\forall x\in\R^d$, $\forall y \in\R$, $\forall z \in\R
^{1\times d}$,
\[
\vert f(t,x,y,z)-f(t',x,y,z)\vert\leq K_{f,t}
\vert t-t'\vert^{1/2}.
\]
Moreover, we set $\eps=Tn^{-a}$ and $N=n^b$, with $a, b \in\R^{+,*}$
two parameters.
Before giving our error estimates, we recall two technical lemmas that
we will prove in Appendices \ref{appA} and \ref{appB}.
\begin{lem}
\label{lemme technique 1}
For all constant $M>0$ there exists a constant $C$ that depends only on
$T$, $M$ and $a$, such that
\[
\prod_{i=0}^{2n-1} (1+Mh_i) \leq C\qquad \forall n \in\N^*.
\]
\end{lem}
\begin{lem}
\label{lemme technique 2}
For all constants $M_1>0$ and $M_2>0$ there exists a constant $C$ that
depends only on $T$, $M_1$, $M_2$ and $a$, such that
\[
\prod_{i=0}^{n-1}\biggl(1+M_1h_i+M_2\frac{h_i}{T-t_{i+1}}\biggr)
\leq Cn^{aM_2}.
\]
\end{lem}

First, we give a convergence result for the Euler scheme.
\begin{prop}
\label{convergence_Euler}
Assume \textup{(HX0)} holds. Then there exists a constant $C$ that does not
depend on $n$, such that
\[
\sup_{0 \leq k \leq2n} \E[ \vert X_{t_k}-X^n_{t_k}
\vert^2
] \leq C\frac{\ln n}{n}.
\]
\end{prop}
\begin{pf}
We just have to copy the classical proof to obtain, thanks to
Lem\-ma~\ref
{lemme technique 1},
\[
\sup_{0 \leq k \leq2n} \E[ \vert X_{t_k}-X^n_{t_k}
\vert^2
] \leq C\sup_{0 \leq i \leq2n-1} h_i=Ch_0.
\]
But
\[
h_0=T(1-n^{-a/n})\leq C\frac{\ln n}{n},
\]
because $(1-n^{-a/n})\sim aT\frac{\ln n}{n}$ when $n \rightarrow
+\infty$,
so the proof is complete.
\end{pf}

Now, let us treat the BSDE approximation. In a first time we will study
the time approximation error on $[T-\eps, T]$.

\begin{prop}
\label{erreur e3 sur [T-eps,T]}
Assume that \textup{(HX0)}, \textup{(HY0)} and \textup{(HY1)} hold. Then there exists
a constant $C$ that does not depend on $n$ and such that
\[
\sup_{n \leq k \leq2n}\E[\vert Y_{t_k}^{N,\eps
,n}-Y_{t_k}^{N,\eps}\vert^2]+ \sum_{k=n}^{2n-1}\E
\biggl[ \int
_{t_k}^{t_{k+1}} \vert Z_{t_k}^{N,\eps,n}-Z_t^{N,\eps}
\vert^2\,dt \biggr]
\leq\frac{C\ln n}{n^{1-2b}}.
\]
%
% (N^2 \E[ \abs{X_{t^X_n}-X^n_{t^X_n}}^2 ]+\frac{
\end{prop}
\begin{pf}
The BSDE (\ref{EDSR_modifiee}) has a linear growth with respect to $z$
on $[T-\eps,T]$ so we are allowed to apply classical results which give
us that
\begin{eqnarray*}
&&\sup_{n \leq k \leq2n}\E[\vert Y_{t_k}^{N,\eps
,n}-Y_{t_k}^{N,\eps}\vert^2]+ \sum_{k=n}^{2n-1}\E
\biggl[ \int
_{t_k}^{t_{k+1}} \vert Z_{t_k}^{N,\eps,n}-Z_t^{N,\eps}
\vert^2\,dt \biggr]\\
&&\qquad \leq C \biggl(\E[ \vert g_N(X_{T})-g_N(X^n_{T})
\vert^2
]+\frac{\eps}{n}\biggr)
\end{eqnarray*}
by using the fact that $g_N$ is $N$-Lipschitz and by applying
Proposition \ref{convergence_Euler}.
\end{pf}
\begin{rem}
\label{remarque discretisation [T-eps,T]}
\begin{longlist}[(1)]
\item[(1)]
When $a \geq1-2b$, then $\eps=Tn^{-a} = o(n^{2b-1}\ln n)$.
We do not need to have a discretization grid on $[T-\eps, T]$; $n+2$
points of discretization are sufficient on $[0,T]$.
\item[(2)] When $a <1-2b$, then it is possible to take only $\lceil
n^c\rceil
$ discretization points on $[T-\eps,T]$ with $a+c=1-2b$. In this case
the error bound becomes
\begin{eqnarray*}
&&\sup_{n \leq k \leq2n}\E[\vert Y_{t_k}^{N,\eps
,n}-Y_{t_k}^{N,\eps}\vert^2]+ \sum_{k=n}^{2n-1}\E
\biggl[ \int
_{t_k}^{t_{k+1}} \vert Z_{t_k}^{N,\eps,n}-Z_t^{N,\eps}
\vert^2\,dt \biggr]\\
&&\qquad\leq C \biggl(\frac{\ln n}{n^{1-2b}}+\frac{1}{n^{a+c}}\biggr)
\end{eqnarray*}
and the Proposition \ref{erreur e3 sur [T-eps,T]} stays true.
\end{longlist}
\end{rem}

Now, let us see what happens on $[0,T-\eps]$.
\begin{theorem}
\label{erreur e3 sur [0,T-eps]}
Assume that \textup{(HX0)}, \textup{(HY0)}, \textup{(HY1)} and \textup{(HX1)} or
\textup{(HX1$'$)} hold. Then for all $\eta>0$, there exists a constant $C$
that does not depend on $N$, $\eps$ and $n$, such that
\[
\sup_{0 \leq k \leq2n}\E[\vert Y_{t_k}^{N,\eps
,n}-Y_{t_k}^{N,\eps}\vert^2]+ \sum_{k=0}^{2n-1}\E
\biggl[ \int
_{t_k}^{t_{k+1}} \vert Z_{t_k}^{N,\eps,n}-Z_t^{N,\eps}
\vert^2\,dt \biggr]
\leq\frac{C}{n^{1-2b-Ka}}
\]
with $K=4(1+\eta)L_{f,z}^2M_{z,2}^2$.
\end{theorem}
\begin{pf}
First, we will study the error on $Y$. From (\ref{EDSR_modifiee}) and
(\ref{Y discretisation temp}) we get
\begin{eqnarray*}
&&Y^{N,\eps}_{t_k}-Y^{N,\eps,n}_{t_k} \\
&&\qquad= \E_{t_k}[Y^{N,\eps
}_{t_{k+1}}-Y^{N,\eps,n}_{t_{k+1}}]\\
&&\qquad\quad{}+ \E_{t_k} \int_{t_k}^{t_{k+1}} \bigl( f(s,X_s,Y_s^{N,\eps
},Z_s^{N,\eps})-f(t_k,X^n_{t_k}, Y^{N,\eps,n}_{t_{k+1}},Z^{N,\eps
,n}_{t_k})\bigr)\,ds.
\end{eqnarray*}
We introduce a parameter $\gamma_k>0$ that will be chosen later. Thanks
to Proposition~\ref{majoration Z modifie} and assumption (HY0), $f$ is
Lipschitz on $[t_k,t_{k+1}]$ with a Lipschitz constant $K_k:=K^1+\frac
{K^2}{(T-t_{k+1})^{1/2}}$ where $K^2=2L_{f,z}M_{z,2} $. A combination
of Young's inequality $(a+b)^2 \leq(1+\gamma_k h_k)a^2+(1+\frac
{1}{\gamma_k h_k})b^2$ and properties of $f$ gives
%
%e16 ###
%
\begin{eqnarray}
\label{preuve convergence 1}
& &\E\vert Y^{N,\eps}_{t_k}-Y^{N,\eps,n}_{t_k}\vert^2
\nonumber\\
&&\qquad\leq (1+\gamma_k h_k) \E\vert\E_{t_k}[Y^{N,\eps
}_{t_{k+1}}-Y^{N,\eps,n}_{t_{k+1}}]\vert^2 \nonumber\\
&&\qquad\quad{}+(1+\eta)^{1/3}K_k^2\biggl(h_k+\frac{1}{\gamma_k}\biggr)\E\int
_{t_k}^{t_{k+1}}\vert Z_s^{N,\eps}-Z_{t_k}^{N,\eps,n}
\vert^2\,ds \\
&&\qquad\quad{} +C\biggl(h_k+\frac{1}{\gamma_k}\biggr)\biggl(h_k^2+\int_{t_k}^{t_{k+1}} \E
\vert X_s-X^n_{t_k}\vert^2\,ds\biggr) \nonumber\\
&&\qquad\quad{}+C\biggl(h_k+\frac{1}{\gamma_k}\biggr)\biggl(\int_{t_k}^{t_{k+1}} \E
\vert Y_s^{N,\eps}-Y^{N,\eps,n}_{t_{k+1}}\vert^2\,ds
\biggr).\nonumber
\end{eqnarray}
We define
\[
\tilde{Z}_{t_k}^{N,\eps,n}:= \frac{1}{h_k}\E_{t_k}[Y^{N,\eps
,n}_{t_{k+1}} (W_{t_{k+1}}-W_{t_k})].
\]
So, $Z_{t_k}^{N,\eps,n}=\rho_{t_{k+1}}(\tilde{Z}_{t_k}^{N,\eps,n})$.
Moreover, Proposition \ref{majoration Z modifie} implies that
$Z_{s}^{N,\eps}=\rho_{t_{k+1}}(Z_{s}^{N,\eps})$ and, since $\rho
_{t_{k+1}}$ is 1-Lipschitz, we have
%
%e17 ###
%
\begin{equation}
\label{rho 1 lipschitz}\quad
\vert Z_s^{N,\eps}-Z_{t_k}^{N,\eps,n}\vert^2 =
\vert\rho_{t_{k+1}}(Z_s^{N,\eps})-\rho_{t_{k+1}}(\tilde
{Z}_{t_k}^{N,\eps,n})\vert^2\leq\vert Z_s^{N,\eps
}-\tilde{Z}_{t_k}^{N,\eps,n}\vert^2.
\end{equation}
As in Theorem \ref{thm de Zhang}, we define $\bar{Z}^{N,\eps}_{t_k}$ by
\begin{eqnarray*}
h_k \bar{Z}^{N,\eps}_{t_k}:\!&=&\E_{t_k} \int
_{t_k}^{t_{k+1}}Z^{N,\eps}_s
\,ds\\
& =& \E_{t_k} \biggl( \biggl(Y^{N,\eps}_{t_{k+1}}+\int_{t_k}^{t_{k+1}}
f(s,X_s,Y_s^{N,\eps},Z_s^{N,\eps})\,ds\biggr)\supmathitt{(W_{t_{k+1}}-W_{t_k})}\biggr).
\end{eqnarray*}
Clearly,
%
%e18 ###
%
\begin{eqnarray}\label{preuve convergence 2}
&&\E\int_{t_k}^{t_{k+1}}\vert Z_s^{N,\eps}-\tilde
{Z}_{t_k}^{N,\eps,n}\vert^2\,ds \nonumber\\[-8pt]\\[-8pt]
&&\qquad= \E\int
_{t_k}^{t_{k+1}}\vert Z_s^{N,\eps}-\bar{Z}_{t_k}^{N,\eps
}\vert^2\,ds
+h_k \E\vert\bar{Z}_{t_k}^{N,\eps}-\tilde
{Z}_{t_k}^{N,\eps,n}\vert^2.\nonumber
\end{eqnarray}
The Cauchy--Schwarz inequality yields
\begin{eqnarray*}
&&\bigl\vert\E_{t_k} \bigl( (Y^{N,\eps}_{t_{k+1}}-Y^{N,\eps
,n}_{t_{k+1}})\supmathitt{(W_{t_{k+1}}-W_{t_{k}})} \bigr)
\bigr\vert^2\\
&&\qquad\leq h_k\{\E_{t_k}(\vert Y^{N,\eps}_{t_{k+1}}-Y^{N,\eps
,n}_{t_{k+1}}\vert^2) - \vert\E_{t_k}(Y^{N,\eps
}_{t_{k+1}}-Y^{N,\eps,n}_{t_{k+1}})\vert^2\}
\end{eqnarray*}
and consequently
%
%e19 ###
%
\begin{eqnarray}\label{preuve convergence 3}
&& h_k \E\vert\bar{Z}_{t_k}^{N,\eps}-\tilde{Z}_{t_k}^{N,\eps
,n}\vert^2
\nonumber\\
&&\qquad \leq (1+\eta)^{1/3} \E[ \E_{t_k}(\vert Y^{N,\eps
}_{t_{k+1}}-Y^{N,\eps,n}_{t_{k+1}}\vert^2) - \vert\E
_{t_k}(Y^{N,\eps}_{t_{k+1}}-Y^{N,\eps,n}_{t_{k+1}})\vert^2
] \\
&&\qquad\quad{} +Ch_k\E\int_{t_k}^{t_{k+1}} \vert f(s,X_s,Y_s^{N,\eps
},Z_s^{N,\eps})\vert^2\,ds.\nonumber
\end{eqnarray}
Plugging (\ref{preuve convergence 2}) and (\ref{preuve convergence 3})
into (\ref{preuve convergence 1}), we get
\begin{eqnarray*}
& & \E\vert Y^{N,\eps}_{t_k}-Y^{N,\eps,n}_{t_k}\vert
^2\\[-2pt]
&&\qquad\leq (1+\gamma_k h_k) \E\vert\E_{t_k}[Y^{N,\eps
}_{t_{k+1}}-Y^{N,\eps,n}_{t_{k+1}}]\vert^2\\[-2pt]
& &\qquad\quad{} +(1+\eta)K_k^2\biggl(h_k+\frac{1}{\gamma_k}\biggr)\E\int
_{t_k}^{t_{k+1}}\vert Z_s^{N,\eps}-\bar{Z}_{t_k}^{N,\eps
}\vert^2\,ds \\[-2pt]
& &\qquad\quad{} +C\biggl(h_k+\frac{1}{\gamma_k}\biggr)\biggl(h_k^2+\int_{t_k}^{t_{k+1}} \E
\vert X_s-X^n_{t_k}\vert^2\,ds\\[-2pt]
&&\qquad\quad\hspace*{74.7pt}{}+\int_{t_k}^{t_{k+1}} \E
\vert Y_s^{N,\eps}-Y^{N,\eps,n}_{t_{k+1}}\vert
^2\,ds\biggr)\\[-2pt]
& &\qquad\quad{} +(1+\eta)^{2/3}K_k^2\biggl(h_k+\frac{1}{\gamma_k}\biggr) \E[ \E
_{t_k}(\vert Y^{N,\eps}_{t_{k+1}}-Y^{N,\eps,n}_{t_{k+1}}
\vert^2)\\[-2pt]
&&\qquad\quad\hspace*{130.9pt}{} - \vert\E_{t_k}(Y^{N,\eps}_{t_{k+1}}-Y^{N,\eps
,n}_{t_{k+1}})\vert^2 ]\\[-2pt]
& &\qquad\quad{} +CK_k^2\biggl(h_k+\frac{1}{\gamma_k}\biggr) h_k\E\int_{t_k}^{t_{k+1}}
\vert f(s,X_s,Y_s^{N,\eps},Z_s^{N,\eps})\vert^2\,ds.
\end{eqnarray*}
Now write
%
%e21 ###
%e20 ###
%
\begin{eqnarray}
\label{Ys-Ynt}
\E\vert Y_s^{N,\eps}-Y^{N,\eps,n}_{t_{k+1}}\vert^2
&\leq&2 \E\vert Y_s^{N,\eps}-Y^{N,\eps}_{t_{k+1}}\vert
^2+2\E\vert Y_{t_{k+1}}^{N,\eps}-Y^{N,\eps,n}_{t_{k+1}}
\vert^2,
\\
\label{Xs-Xnt}
\E\vert X_s-X^n_{t_k}\vert^2 &\leq&2 \E\vert
X_s-X_{t_k}\vert^2+2 \E\vert X_{t_k}-X^n_{t_k}
\vert^2
\end{eqnarray}
and we obtain
\begin{eqnarray*}
&& \E\vert Y^{N,\eps}_{t_k}-Y^{N,\eps,n}_{t_k}\vert^2\\
&&\qquad\leq (1+\gamma_k h_k) \E\vert\E_{t_k}[Y^{N,\eps
}_{t_{k+1}}-Y^{N,\eps,n}_{t_{k+1}}]\vert^2\\
&&\qquad\quad{} +(1+\eta)K_k^2\biggl(h_k+\frac{1}{\gamma_k}\biggr)\E\int
_{t_k}^{t_{k+1}}\vert Z_s^{N,\eps}-\bar{Z}_{t_k}^{N,\eps
}\vert^2\,ds \\
&&\qquad\quad{} +C\biggl(h_k+\frac{1}{\gamma_k}\biggr)\biggl(h_k^2+\int_{t_k}^{t_{k+1}} \E
\vert X_s-X_{t_k}\vert^2\,ds+h_k\E\vert
X_{t_k}-X^n_{t_k}\vert^2 \biggr)\\
& &\qquad\quad{} +C\biggl(h_k+\frac{1}{\gamma_k}\biggr)\biggl( \int_{t_k}^{t_{k+1}} \E
\vert Y_s^{N,\eps}-Y^{N,\eps}_{t_{k+1}}\vert^2\,ds+h_k \E
\vert Y_{t_{k+1}}^{N,\eps}-Y^{N,\eps,n}_{t_{k+1}}\vert^2
\biggr)\\
& &\qquad\quad{} +(1+\eta)^{2/3}K_k^2\biggl(h_k+\frac{1}{\gamma_k}\biggr) \E[ \E
_{t_k}(\vert Y^{N,\eps}_{t_{k+1}}-Y^{N,\eps,n}_{t_{k+1}}
\vert^2)\\
&&\qquad\quad\hspace*{131pt}{} - \vert\E_{t_k}(Y^{N,\eps}_{t_{k+1}}-Y^{N,\eps
,n}_{t_{k+1}})\vert^2 ]\\
& &\qquad\quad{} +CK_k^2\biggl(h_k+\frac{1}{\gamma_k}\biggr) h_k\E\int_{t_k}^{t_{k+1}}
\vert f(s,X_s,Y_s^{N,\eps},Z_s^{N,\eps})\vert^2\,ds.
\end{eqnarray*}
Taking $\gamma_k=(1+\eta)^{2/3}K_k^2$ and for $h_k$ small enough, it gives
\begin{eqnarray*}
&& \E\vert Y^{N,\eps}_{t_k}-Y^{N,\eps,n}_{t_k}\vert^2\\[-2pt]
&&\qquad\leq \bigl(1+Ch_ k+(1+\eta)^{2/3}K_k^2h_k\bigr)\E\vert Y^{N,\eps
}_{t_{k+1}}-Y^{N,\eps,n}_{t_{k+1}}\vert^2+Ch_k^2\\[-2pt]
&&\qquad\quad{}+Ch_k\max_{0\leq k\leq n} \E\vert X_{t_k}-X^n_{t_k}
\vert^2 \\[-2pt]
& &\qquad\quad{} +C \E\int_{t_k}^{t_{k+1}}\vert Z_s^{N,\eps}-\bar
{Z}_{t_k}^{N,\eps}\vert^2\,ds+ C \int_{t_k}^{t_{k+1}} \E
\vert X_s-X_{t_k}\vert^2\,ds\\[-2pt]
& &\qquad\quad{} +C \int_{t_k}^{t_{k+1}} \E\vert Y_s^{N,\eps}-Y^{N,\eps
}_{t_{k+1}}\vert^2\,ds +Ch_k\E\int_{t_k}^{t_{k+1}}
f(s,X_s,Y_s^{N,\eps
},Z_s^{N,\eps})^2\,ds,
\end{eqnarray*}
because $K_k^2h_k \leq C(h_0+h_k(T-t_{k+1})^{-1}) \leq C\frac
{\ln n}{n}$. The Gronwall's lemma gives us
\begin{eqnarray*}
\hspace*{-4pt}&& \E\vert Y^{N,\eps}_{t_k}-Y^{N,\eps,n}_{t_k}\vert^2\\[-2pt]
\hspace*{-4pt}&&\qquad\leq C\sum_{j=0}^{n-1}\Biggl[\prod_{i=0}^{j-1}\bigl(1+Ch_ i+(1+\eta
)^{2/3}K_i^2h_i\bigr) \Biggr]\\[-2pt]
\hspace*{-4pt}&&\qquad\quad\hspace*{23pt}{}\times\biggl[h_j^2+h_j\max_{0\leq l\leq n}
\E\vert X_{t_l}-X^n_{t_l}\vert^2 \\[-2pt]
\hspace*{-4pt}&&\qquad\quad\hspace*{39.3pt}{} + \E\int_{t_j}^{t_{j+1}}(\vert Z_s^{N,\eps}-\bar
{Z}_{t_j}^{N,\eps}\vert^2 + \vert X_s-X_{t_j}
\vert^2+\vert Y_s^{N,\eps}-Y^{N,\eps}_{t_{j+1}}\vert
^2) \,ds\\[-2pt]
\hspace*{-4pt}& &\qquad\quad\hspace*{145.3pt}{} +h_j\E\int_{t_j}^{t_{j+1}} \vert
f(s,X_s,Y_s^{N,\eps},Z_s^{N,\eps})\vert^2\,ds\biggr]\\[-2pt]
\hspace*{-4pt}& &\qquad\quad{} + \Biggl[\prod_{i=0}^{n-1}\bigl(1+Ch_ i+(1+\eta)^{2/3}K_i^2h_i\bigr)
\Biggr]\E
\vert Y^{N,\eps}_{t_n}-Y^{N,\eps,n}_{t_n}\vert^2.
\end{eqnarray*}

Then, we apply Lemma \ref{lemme technique 2}.
\begin{eqnarray*}
\hspace*{-4pt}&& \E\vert Y^{N,\eps}_{t_k}-Y^{N,\eps,n}_{t_k}\vert^2\\[-2pt]
\hspace*{-4pt}&&\qquad\leq Cn^{(1+\eta)(K^2)^2a}\\[-2pt]
\hspace*{-4pt}&&\qquad\quad{}\times\Biggl[h_0+\max_{0\leq l\leq
n} \E\vert X_{t_l}-X^n_{t_l}\vert^2 \\[-2pt]
\hspace*{-4pt}& &\hspace*{15.64pt}\qquad\quad{} + \sum_{j=0}^n \E\biggl(\int_{t_j}^{t_{j+1}}\vert
Z_s^{N,\eps}-\bar{Z}_{t_j}^{N,\eps}\vert^2 + \vert
X_s-X_{t_j}\vert^2+\vert Y_s^{N,\eps}-Y^{N,\eps
}_{t_{j+1}}\vert^2 \,ds\biggr)\\[-2pt]
\hspace*{-4pt}& &\hspace*{55.1pt}\qquad\quad{} +h_0\E\int_{0}^{t_{n}} \vert f(s,X_s,Y_s^{N,\eps
},Z_s^{N,\eps})\vert^2\,ds + \E\vert Y^{N,\eps
}_{t_n}-Y^{N,\eps,n}_{t_n}\vert^2\Biggr].
\end{eqnarray*}

A classical estimation gives us $\E[\vert
X_s-X_{t_j}\vert^2
]\leq\vert s-t_j\vert$. Moreover, since $Z^{N,\eps}$ is bounded,
\begin{eqnarray*}
&&\E\int_{0}^{t_{n}} \vert f(s,X_s,Y_s^{N,\eps},Z_s^{N,\eps
})\vert^2\,ds\\
&&\qquad\leq CT(1+\vert Y^{N,\eps}\vert_{\infty})+C\E
\biggl[\int_{0}^{t_{n}}
\vert Z_s^{N,\eps}\vert^4\,ds\biggr]\\
&&\qquad \leq CT(1+\vert Y^{N,\eps}\vert_{\infty})+Cn^{2b}\E
\biggl[\int
_{0}^{T} \vert Z_s^{N,\eps}\vert^2\,ds\biggr].
\end{eqnarray*}
But we have an a priori estimate for $\E[\int_{0}^{T}
\vert Z_s^{N,\eps}\vert^2\,ds]$ that does not depend on
$N$ and $\eps$. So
%
%e22 ###
%
\begin{equation}
\label{majoration int fcarre}
\E\int_{0}^{t_{n}} \vert f(s,X_s,Y_s^{N,\eps},Z_s^{N,\eps
})\vert^2\,ds
\leq Cn^{2b}.
\end{equation}
With the same type of argument we also have
%
%e23 ###
%
\begin{equation}
\label{Ys-Yt}
\E\vert Y_s^{N,\eps}-Y^{N,\eps}_{t_{j+1}}\vert^2 \leq
Ch_jn^{2b}.
\end{equation}
If we add Zhang's path regularity Theorem \ref{thm de Zhang},
Propositions \ref{convergence_Euler} and \ref{erreur e3 sur
[T-eps,T]}, we finally obtain
%
%e24 ###
%
\begin{equation}
\label{erreur provisoire sur Y}
\E\vert Y^{N,\eps}_{t_k}-Y^{N,\eps,n}_{t_k}\vert^2
\leq Cn^{(1+\eta
)(K^2)^2a}\frac{n^{2b}\ln n}{n} = C\frac{\ln n}{n^{1-2b-(1+\eta)(K^2)^2a}}.
\end{equation}
Now, let us deal with the error on $Z$. First of all, (\ref{rho 1
lipschitz}) gives us
\[
\sum_{k=0}^{n-1}\E\biggl[ \int_{t_k}^{t_{k+1}} \vert
Z_{t_k}^{N,\eps,n}-Z_t^{N,\eps}\vert^2\,dt \biggr] \leq\sum
_{k=0}^{n-1}\E\biggl[ \int
_{t_k}^{t_{k+1}} \vert\tilde{Z}_{t_k}^{N,\eps,n}-Z_t^{N,\eps
}\vert^2\,dt
\biggr].
\]
For $0\leq k\leq n-1$, we can use (\ref{preuve convergence
2}) and (\ref{preuve convergence 3}) to obtain
\begin{eqnarray*}
&& \E\biggl[ \int_{t_k}^{t_{k+1}} \vert\tilde
{Z}_{t_k}^{N,\eps,n}-Z_t^{N,\eps}\vert^2\,dt\biggr]\\
&&\qquad\leq \E\biggl[ \int_{t_k}^{t_{k+1}} \vert\bar
{Z}_{t_k}^{N,\eps}-Z_t^{N,\eps}\vert^2\,dt \biggr]\\
&&\qquad\quad{}+(1+\eta)^{2/3} \E[ \E_{t_k}(\vert Y^{N,\eps
}_{t_{k+1}}-Y^{N,\eps,n}_{t_{k+1}}\vert^2) - \vert\E
_{t_k}(Y^{N,\eps}_{t_{k+1}}-Y^{N,\eps,n}_{t_{k+1}})\vert^2
]\\
&&\qquad\quad{}+Ch_k\E\biggl[\int_{t_k}^{t_{k+1}} \vert f(s,X_s,Y_s^{N,\eps
},Z_s^{N,\eps})\vert^2\,ds\biggr].
\end{eqnarray*}
Inequality (\ref{majoration int fcarre}) and estimates for $Z$ give us
%
%e25 ###
%
\begin{eqnarray}\label{preuve convergence 4}
&& \sum_{k=0}^{n-1}\E\biggl[ \int_{t_k}^{t_{k+1}} \vert
Z_{t_k}^{N,\eps,n}-Z_t^{N,\eps}\vert^2\,dt \biggr]\nonumber\\
&&\qquad\leq \sum_{k=0}^{n-1}\E\biggl[ \int_{t_k}^{t_{k+1}} \vert
\bar{Z}_{t_k}^{N,\eps}-Z_t^{N,\eps}\vert^2\,dt
\biggr]\nonumber\\
&&\qquad\quad{}+(1+\eta)^{2/3} \sum_{k=0}^{n-1}\E[ \E_{t_k}(\vert
Y^{N,\eps}_{t_{k+1}}-Y^{N,\eps,n}_{t_{k+1}}\vert^2) -
\vert\E_{t_k}(Y^{N,\eps}_{t_{k+1}}-Y^{N,\eps,n}_{t_{k+1}})
\vert^2 ]\nonumber\\
&&\qquad\quad{}+ Ch_0\E\biggl[\int_0^T\vert f(s,X_s,Y_s^{N,\eps
},Z_s^{N,\eps})\vert^2\,ds\biggr]\\
&&\qquad\leq \sum_{k=0}^{n-1}\E\biggl[ \int_{t_k}^{t_{k+1}} \vert
\bar{Z}_{t_k}^{N,\eps}-Z_t^{N,\eps}\vert^2\,dt
\biggr]\nonumber\\
&&\qquad\quad{}+(1+\eta)^{2/3} \sum_{k=0}^{n-1}\E[ \E_{t_k}(\vert
Y^{N,\eps}_{t_{k}}-Y^{N,\eps,n}_{t_{k}}\vert^2) -
\vert\E_{t_k}(Y^{N,\eps}_{t_{k+1}}-Y^{N,\eps,n}_{t_{k+1}})
\vert^2 ]\nonumber\\
&&\qquad\quad{}+C\E[\vert Y^{N,\eps}_{t_{n}}-Y^{N,\eps
,n}_{t_{n}}\vert^2]+
Ch_0n^{2b}\nonumber
\end{eqnarray}
with an index change in the penultimate line. Then, by using (\ref
{preuve convergence 1}) we get
%
%e26 ###
%
\begin{eqnarray}
\label{preuve convergence 5}
& & (1+\eta)^{2/3} \E[ \E_{t_k}(\vert Y^{N,\eps
}_{t_{k}}-Y^{N,\eps,n}_{t_{k}}\vert^2) - \vert\E
_{t_k}(Y^{N,\eps}_{t_{k+1}}-Y^{N,\eps,n}_{t_{k+1}})\vert^2
]\nonumber\\
&&\qquad\leq C\gamma_k h_k \E\vert\E_{t_k}[Y^{N,\eps
}_{t_{k+1}}-Y^{N,\eps,n}_{t_{k+1}}]\vert^2\nonumber\\[-8pt]\\[-8pt]
&&\qquad\quad{}+(1+\eta)K_k^2\biggl(h_k+\frac{1}{\gamma_k}\biggr)\E\int
_{t_k}^{t_{k+1}}\vert Z_s^{N,\eps}-Z_{t_k}^{N,\eps,n}
\vert^2\,ds \nonumber\\
&&\qquad\quad{}+ C\biggl(h_k+\frac{1}{\gamma_k}\biggr)h_k\Bigl(h_k+\sup_{s \in
[t_k,t_{k+1}]} \E
[\vert X_s-X^n_{t_k}\vert^2+ \vert
Y_s^{N,\eps}-Y^{N,\eps,n}_{t_{k+1}}\vert^2]
\Bigr).\nonumber
\end{eqnarray}
Thanks to (\ref{Ys-Ynt}), (\ref{Xs-Xnt}), (\ref{Ys-Yt}) and a classical
estimation on $\E[\vert X_s-X_{t_k}\vert^2
]$ we have
\begin{eqnarray*}
&&\sup_{s \in[t_k,t_{k+1}]} \E[\vert X_s-X^n_{t_k}
\vert^2+ \vert Y_s^{N,\eps}-Y^{N,\eps,n}_{t_{k+1}}\vert
^2]\\
&&\qquad\leq C
(h_kn^{2b}+ \E[\vert Y_{t_{k+1}}^{N,\eps}-Y^{N,\eps
,n}_{t_{k+1}}\vert^2]).
\end{eqnarray*}
Let us set $\gamma_k=3(1+\eta)K_k^2$. We recall that $h_kK_k^2
\leq\frac{C\ln n}{n} \rightarrow0$ when $n\rightarrow0$. So,
for $n$ big enough, (\ref{preuve convergence 5}) becomes
\begin{eqnarray*}
& & (1+\eta)^{2/3} \E[ \E_{t_k}(\vert Y^{N,\eps
}_{t_{k}}-Y^{N,\eps,n}_{t_{k}}\vert^2) - \vert\E
_{t_k}(Y^{N,\eps}_{t_{k+1}}-Y^{N,\eps,n}_{t_{k+1}})\vert^2
]\\
&&\qquad\leq \frac{C\ln n}{n} \E[\vert Y^{N,\eps
}_{t_{k+1}}-Y^{N,\eps,n}_{t_{k+1}}\vert^2]+\frac
{1}{2}\E\int
_{t_k}^{t_{k+1}}\vert Z_s^{N,\eps}-Z_{t_k}^{N,\eps,n}
\vert^2\,ds\\
&&\qquad\quad{}+ Ch_0h_kn^{2b}.
\end{eqnarray*}
If we inject this last estimate in (\ref{preuve convergence 4}) and we
use Theorem \ref{thm de Zhang}, we obtain
\begin{eqnarray*}
&& \frac{1}{2}\sum_{k=0}^{n-1}\E\biggl[ \int_{t_k}^{t_{k+1}}
\vert Z_{t_k}^{N,\eps,n}-Z_t^{N,\eps}\vert^2\,dt \biggr]\\
&&\qquad\leq
Ch_0n^{2b}+C\ln n \sup_{0 \leq k \leq n-1} \E[\vert
Y^{N,\eps}_{t_{k+1}}-Y^{N,\eps,n}_{t_{k+1}}\vert^2].
\end{eqnarray*}
By using (\ref{erreur provisoire sur Y}) and Proposition \ref{erreur e3
sur [T-eps,T]}, we finally have
\begin{eqnarray*}
&&\sup_{0 \leq k \leq2n}\E[\vert Y_{t_k}^{N,\eps
,n}-Y_{t_k}^{N,\eps}\vert^2]+ \sum_{k=0}^{2n-1}\E
\biggl[ \int
_{t_k}^{t_{k+1}} \vert Z_{t_k}^{N,\eps,n}-Z_t^{N,\eps}
\vert^2\,dt \biggr]\\
&&\qquad\leq C\frac{(\ln n)^2}{n^{1-2b-Ka}}
\end{eqnarray*}
with\vspace*{1pt} $K=4(1+\eta)L_{f,z}^2M_{z,2}^2$. Since this estimate is true for
every $\eta>0$, we have proved the result.
\end{pf}

%s4.3 ###
\subsection{Study of the global error $e(N,\eps,n)$}\label{sec4.3}
Let us study errors $e_1(N)$ and $e_2(N,\eps)$.
\begin{prop}
\label{erreur approximation eps}
Let us assume that \textup{(HX0)} and \textup{(HY0)} hold. There exists a
constant $C>0$ such that
\[
e_2(N,\eps) \leq\frac{C}{n^{2a-4b}}.
\]
\end{prop}
\begin{pf}
We just have to notice that
\[
\vert f(t,X_t,Y_t^{N,\eps},Z_t^{N,\eps})-f(t,X_t,Y_t^{N,\eps
},0)\vert\leq C\vert Z_t^{N,\eps}\vert^2
\]
and $\vert Z_t^{N,\eps}\vert$ is bounded by $Cn^b$.
\end{pf}

For $g_N$ we use the classical Lipschitz approximation
\[
g_N(x)=\inf\{g(u)+N\vert x-u\vert| u \in\R
^d\}.
\]

\begin{prop}
\label{erreur approximation g alpha holder}
We assume that \textup{(HX0)} holds and $g$ is $\alpha$-H\"{o}lder. Then,
there exists a constant $C$ such that
\[
e_1(N) \leq\frac{C}{n^{{2b\alpha}/({1-\alpha})}}.
\]
\end{prop}
\begin{pf}
$g_N$ is a $N$-Lipschitz function and $g_N \rightarrow g$ when $N
\rightarrow+ \infty$ uniformly on $\R^d$. More precisely, we have
\[
\vert g-g_N\vert_{\infty}\leq\frac{C}{N^{{\alpha
}/({1-\alpha})}}.
\]
\upqed\end{pf}
\begin{rem}
\label{exemple x^alpha}
For some explicit examples, it is possible to have a better
convergence speed. For example, let us take $g(x)=(\vert x
\vert^{\alpha}
\mathbh{1}_{x\geq0})\wedge C$ and assume that $\sigma$ is
invertible. Then, we can use the fact that this function is not
Lipschitz only in $0$ and obtain
\[
e_1(N) \leq\frac{C}{n^{{2\alpha b}/({1-\alpha})}} {\mathbb
P}\bigl(X_T
\in\bigl[0, N^{{-1}/({1-\alpha})}\bigr]\bigr)^{1/q}\leq\frac
{C}{n^{({b}/({1-\alpha}))(2\alpha+{1}/{q})}}.
\]
\end{rem}
\begin{rem}
\label{exemple_indicatrice}
It is also possible to obtain convergence speed when $g$ is not
$\alpha$-H\"{o}lder. For example, we assume that $\sigma$ is invertible
and we set $g(x)=\prod_{i=1}^d \mathbh{1}_{x_i>0}(x)$. Then
\[
e_1(N) \leq C\Biggl[\sum_{i=1}^{d} {\mathbb P}\bigl((X_T)_i \in
[0,1/N]\bigr)
\Biggr]^{1/q} \leq\frac{C}{N^{1/q}}=\frac{C}{n^{b/q}}.
\]
\end{rem}

Now we are able to gather all these errors.
\begin{theorem}
\label{erreur_globale}
We assume that \textup{(HX0)}, \textup{(HY0)}, \textup{(HY1)} and \textup{(HX1)} or
\textup{(HX1$'$)} hold. We assume also that $g$ is $\alpha$-H\"{o}lder. Then
for all $\eta>0$, there exists a constant $C>0$ that does not depend on
$n$ such that
\[
e(n):=e(N,\eps,n)\leq\frac{C}{n^{{2\alpha}/({(2-\alpha
)(2+K)-2+2\alpha})}}
\]
with $K=4(1+\eta)L_{f,z}^2M_{z,2}^2$.
\end{theorem}
\begin{pf}
Thanks to Theorem \ref{erreur e3 sur [0,T-eps]}, Propositions \ref
{erreur approximation eps} and \ref{erreur approximation g
alpha holder} we have
\[
e(n) \leq\frac{C}{n^{1-2b-Ka}}+\frac{C}{n^{2a-4b}}+\frac
{C}{n^{{2\alpha b}/({1-\alpha})}}.\vadjust{\goodbreak}
\]
Then we only need to set $a:=\frac{1+2b}{2+K}$ and $b:=\frac{1-\alpha
}{(2-\alpha)(2+K)-2+2\alpha}$ to obtain the result.
%\rightqed
\end{pf}
\begin{cor}
\label{corollaire f sous-quadratique}
We assume that assumptions of Theorem \ref{erreur_globale} hold.
Moreover, we assume that $f$ has a sub-quadratic growth with respect to
$z$; there exists $0<\beta<1$ such that, for all $t \in[0,T]$, $x\in
\R
^d$, $y \in R$, $z,z' \in\R^{1\times d}$,
\[
\vert f(t,x,y,z)-f(t,x,y,z')\vert\leq
\bigl(K_{f,z}+L_{f,z}(\vert z\vert^{\beta
}+\vert z'\vert^{\beta})\bigr)\vert z-z'\vert.
\]
Then we are allowed to take $K$ as small as we want. So, for all $\eta
>0$, there exists a constant $C>0$ that does not depend on $n$ such that
\[
e(n) \leq\frac{C}{n^{\alpha-\eta}}.
\]
\end{cor}
%
% When we are allowed to take $K$ as small as we want, then we have $
%have a discretization grid on $[T-\eps,T]$: $n+2$ points of
%discretization are sufficient on %$[0,T]$.
%
\begin{rem}
We are able to specify Remark \ref{remarque discretisation
[T-eps,T]} in our case, when $a=\frac{1+2b}{2+K}$ and $b=\frac
{1-\alpha
}{(2-\alpha)(2+K)-2+2\alpha}$.
\begin{longlist}[(1)]
\item[(1)] When $K \leq\frac{2-3\alpha}{2-\alpha}$, that is to say,
when $\alpha< 2/3$ and $K$ is sufficiently small, then we do not need
to have a discretization grid on $[T-\eps,T]$.
\item[(2)] When $K > \frac{2-3\alpha}{2-\alpha}$, then it is possible to
take only $\lceil n^c \rceil$ discretization points on $[T-\eps,T]$ with
\[
c = 1+\frac{3\alpha-4}{(2-\alpha)(2+K)-2+2\alpha}.
\]
\end{longlist}
\end{rem}

Theorem \ref{erreur_globale} is not interesting in practice because the
speed of convergence depends strongly on $K$. But we\vspace*{1pt} see that the
global error becomes $e(n) \leq\frac{C}{n^{\alpha-\eta}}$ when we
are allowed to choose $K$ as small as we want. Under extra assumption
we can show that we are allowed to take the constant $M_{z,2}$ as small
as we want.

\subsubsection*{\textup{(HX2)}}
$b$ is bounded on $[0,T]\times\R^d$ by a constant $M_b$.
\begin{theorem}
\label{vitesse_optimale}
We assume that \textup{(HX0)}, \textup{(HY0)}, \textup{(HY1)}, \textup{(HX2)} and
\textup{(HX1)} or \textup{(HX1$'$)} hold. We assume also that $g$ is $\alpha$-H\"
{o}lder. Then for all $\eta>0$, there exists a constant $C>0$ that does
not depend on $n$ such that
\[
e(n) \leq\frac{C}{n^{\alpha-\eta}}.
\]
\end{theorem}
\begin{rem}
With the assumptions of the previous theorem, it is also possible
to have an estimate of the global error for examples given in
Remarks \ref{exemple x^alpha} and \ref{exemple_indicatrice}. When
$g(x)=(\vert x\vert^{\alpha} \mathbh{1}_{x\geq
0})\wedge C$, we have
\[
e(n) \leq\frac{C}{n^{\alpha+({1-\alpha})/({1+2q})-\eta}}
\]
and when $g(x)=\prod_{i=1}^d \mathbh{1}_{x_i>0}(x)$, we have
\[
e(n) \leq\frac{C}{n^{{1}/({1+2q})-\eta}}.
\]
\end{rem}
\begin{pf*}{Proof of Theorem \ref{vitesse_optimale}}
First, we suppose that $f$ is differentiable with respect to $z$.
Thanks to Remark \ref{estimation temporelle plus precise pour Z} we see
that it is sufficient to show that
\[
\E^{\Q^{N,\eps}}\biggl[\int_t^T \vert Z^{N,\eps}_s
\vert^2 \,ds
\big|\mathcal
{F}_t\biggr]
\]
is small uniformly in $\omega$, $N$ and $\eps$ when $t$ is close to $T$.
We will obtain an estimation for this quantity by applying the same
computation as \cite{Briand-Confortola-08} for the BMO norm estimate of
$Z$, page 831. Thus, we have
\[
\E^{\Q^{N,\eps}}\biggl[\int_t^T \vert Z^{N,\eps}_s
\vert^2 \,ds
\big|\mathcal
{F}_t\biggr] \leq\E^{\Q^{N,\eps}}[ \vert\varphi
(Y^{N,\eps}_T)-\varphi(Y^{N,\eps}_t)\vert|\mathcal
{F}_t]+C(T-t)
\]
with $\varphi(x)=(e^{2c(x+m)}-2c(x+m)-1)/(2c^2)$, $m=\vert
Y\vert_{\infty}$
and $c$ that depends on constants in assumption (HY0) but does not
depend on $\nabla_z f$. Let us notice that $m$, $c$ and so $\varphi$ do
not depend on $N$ and $\eps$. Since $Y$ is bounded, $\varphi$ is a
Lipschitz function, so
\[
\E^{\Q^{N,\eps}}\biggl[\int_t^T \vert Z^{N,\eps}_s
\vert^2 \,ds
\big|\mathcal
{F}_t\biggr] \leq C\E^{\Q^{N,\eps}}[ \vert Y^{N,\eps
}_T-Y^{N,\eps}_t\vert|\mathcal{F}_t]+C(T-t).
\]
We denote by $(Y^{N,\eps,t,x},Z^{N,\eps,t,x})$ the solution of
BSDE (\ref{EDSR_modifiee}) when $X^{t,x}_t=x$. As usual, we set
$X^{t,x}_s=x$ and $Z^{N,\eps,t,x}_s=0$ for $s\leq t$ and we define
$u^{N,\eps}(t,x):=Y_t^{N,\eps,t,x}$. Then we give a proposition that we
will prove in Appendix \ref{appC}.
\begin{prop}
\label{u uniformement continue}
We assume that \textup{(HX0)}, \textup{(HY0)}, \textup{(HY1)}, \textup{(HX2)} and
\textup{(HX1)} or \textup{(HX1$'$)} hold. We assume also that $g$ is uniformly
continuous on $\R^d$. Then $u^{N,\eps}$ is uniformly continuous on
$[0,T] \times\R^d$ and there exists $\omega$ a concave modulus of
continuity for all functions in $\{u^{N,\eps}| N \in\N, \eps
>0\}$,
that is, $\omega$ does not depend on $N$ and $\eps$.
\end{prop}

Then
\begin{eqnarray*}
&& \E^{\Q^{N,\eps}}[\vert Y^{N,\eps}_T-Y^{N,\eps
}_t\vert
|\mathcal
{F}_t]\\
&&\qquad= \E^{\Q^{N,\eps}}[\vert u^{N,\eps}(T,X_T)-u^{N,\eps
}(t,X_t)\vert
|\mathcal{F}_t]\\
&&\qquad \leq \E^{\Q^{N,\eps}}\bigl[ \mathbh{1}_{\vert\int_t^T
\sigma(s)\,d\tilde{W}_s\vert\leq\nu} \vert u^{N,\eps
}(T,X_T)-u^{N,\eps}(t,X_t)\vert\\
& &\qquad\quad\hspace*{84pt}{} + 2\vert Y^{N,\eps}\vert_{\infty}\mathbh
{1}_{\vert\int_t^T \sigma(s)\,d\tilde{W}_s\vert>\nu}
|\mathcal{F}_t\bigr]\\
&&\qquad \leq \E^{\Q^{N,\eps}}\bigl[\omega\bigl(\vert T-t
\vert+\mathbh
{1}_{\vert\int_t^T \sigma(s)\,d\tilde{W}_s\vert\leq\nu
}\vert X_T-X_t\vert\bigr)\\
&&\qquad\quad\hspace*{64.6pt}{} + 2\vert Y^{N,\eps}\vert_{\infty}\mathbh
{1}_{\vert\int_t^T \sigma(s)\,d\tilde{W}_s\vert>\nu}
|\mathcal{F}_t\bigr]
\end{eqnarray*}
with $d\tilde{W}_s = dW_s-\nabla_z f^{\eps}(s,X_s,Y^{N,\eps
}_s,Z^{N,\eps
}_s)\,ds$. But,
\begin{eqnarray*}
& &\mathbh{1}_{\vert\int_t^T \sigma(s)\,d\tilde{W}_s
\vert\leq\nu}\vert X_T-X_t\vert\\
&&\qquad= \mathbh{1}_{\vert\int_t^T \sigma(s)\,d\tilde{W}_s
\vert\leq\nu}
\biggl|\int_t^T b(s,X_s)\,ds\\
&&\hspace*{67.8pt}\qquad\quad{}+\int_t^T \nabla_z f^{\eps}(s,X_s,Y^{N,\eps}_s,Z^{N,\eps}_s)\,ds+
\int_t^T \sigma(s)\,d\tilde{W}_s\biggr|\\
&&\qquad\leq M_b(T-t)+\nu+ C\int_t^T (1+\vert Z^{N,\eps}_s
\vert)\,ds\\
&&\qquad\leq C(T-t)+\nu+ C(T-t)^{1/2}\biggl(\int_t^T \vert
Z^{N,\eps}_s\vert^2\,ds\biggr)^{1/2}.
\end{eqnarray*}
Since $\omega$ is concave, we have by Jensen's inequality
\begin{eqnarray*}
& &\E^{\Q^{N,\eps}}\bigl[\omega\bigl(\vert T-t\vert
+\mathbh{1}_{\vert\int_t^T \sigma(s)\,d\tilde{W}_s
\vert\leq\nu}\vert X_T-X_t\vert\bigr)
|\mathcal{F}_t\bigr]\\
&&\qquad\leq \omega\biggl( C\vert T-t\vert+\nu
+C(T-t)^{1/2}\E^{\Q^{N,\eps
}}\biggl[\biggl(\int_t^T \vert Z^{N,\eps}_s\vert
^2\,ds\biggr)^{1/2}
\Big|\mathcal
{F}_t\biggr]\biggr)\\
&&\qquad\leq \omega\biggl( C\vert T-t\vert+\nu
+C(T-t)^{1/2}\E^{\Q^{N,\eps
}}\biggl[\int_t^T \vert Z^{N,\eps}_s\vert^2\,ds
\big|\mathcal{F}_t
\biggr]^{1/2}\biggr)\\
&&\qquad\leq \omega\bigl( C\vert T-t\vert+\nu
+C(T-t)^{1/2}\Vert Z^{N,\eps}\Vert_{\mathrm{BMO}(\Q)}\bigr).
\end{eqnarray*}
But, $\Vert Z^{N,\eps}\Vert_{\mathrm{BMO}(\Q)}$ only depends on
constants in
assumption (HY0), so it is bounded uniformly in $N$ and $\eps$.
Moreover, $\vert\int_t^T \sigma(s)\,d\tilde{W}_s\vert$
is independent of
$\mathcal{F}_t$ so we have by the Markov inequality
\begin{eqnarray*}
\E^{\Q^{N,\eps}}\bigl[\mathbh{1}_{\vert\int_t^T \sigma
(s)\,d\tilde{W}_s\vert>\nu} |\mathcal{F}_t\bigr] &=&
\Q^{N,\eps}\biggl( \biggl\vert\int_t^T \sigma(s)\,d\tilde
{W}_s\biggr\vert>\nu\biggr)\\
%&\leq& \frac{\int_t^T \tr{\sigma}(s)\sigma(s)\,ds}{\nu^2}\\
&\leq& \frac{C(T-t)^{1/2}}{\nu}.
\end{eqnarray*}
Finally, we have
\begin{eqnarray*}
\E^{\Q^{N,\eps}}[\vert Y^{N,\eps}_T-Y^{N,\eps}_t
\vert
|\mathcal
{F}_t]& \leq& \omega(C\vert T-t\vert
^{1/2}+\nu) +
C\frac{(T-t)^{1/2}}{\nu}\\
&\leq& \omega(C\vert T-t\vert^{1/2}+\vert
T-t\vert^{1/4}) +
C\vert T-t\vert^{1/4}
\end{eqnarray*}
by setting $\nu=\vert T-t\vert^{1/4}$ and $\E^{\Q
^{N,\eps}}[\vert Y^{N,\eps}_T-Y^{N,\eps}_t\vert
|\mathcal{F}_t]\rightarrow0$
uniformly in $\omega$, $N$ and $\eps$ when $t\rightarrow T$. When $f$
is not differentiable with respect to $z$ but is only locally
Lipschitz, then we can prove the result by a standard
approximation.\vadjust{\goodbreak}
\end{pf*}

\begin{appendix}
%s5 ###
\section{\texorpdfstring{Proof of Lemma
\protect\lowercase{\ref{lemme technique 1}}}{Proof of Lemma 4.4}}
\label{appA}
We have
\[
\prod_{i=0}^{2n-1} (1+Mh_i) = \Biggl(\prod_{i=0}^{n-1} (1+Mh_i)
\Biggr)\Biggl(\prod_{i=n}^{2n-1} (1+Mh_i)\Biggr).
\]
First,
\[
\prod_{i=n}^{2n-1} (1+Mh_i)\leq\biggl(1+M\frac{T}{n}\biggr)^n
\leq C.
\]
Moreover, for $0 \leq i \leq n-1$,
\[
h_i=t_{i+1}-t_i=Tn^{-ai/n}\bigl(1-e^{-({a\ln n})/{n}}\bigr)\leq Tn^{-ai/n}
a\frac{\ln n}{n}
\]
thanks to the convexity of the exponential function. So
\begin{eqnarray*}
\prod_{i=0}^{n-1} (1+Mh_i) &\leq& \prod_{i=0}^{n-1}
\biggl(1+MTan^{-ai/n}\frac{\ln n}{n}\biggr)\\
& = & \exp\Biggl( \sum_{i=0}^{n-1} \ln\biggl( 1+MTan^{-ai/n}\frac
{\ln
n}{n}\biggr)\Biggr)\\
& \leq& \exp\Biggl( \sum_{i=0}^{n-1} MTa(n^{-a/n})^i\frac{\ln
n}{n}\Biggr)\\
& \leq& \exp\biggl( MTa\frac{\ln n}{n}\biggl(\frac
{1-(1/n^a)}{1-(1/n^{(a/n)})}\biggr)\biggr)\\
& \leq& \exp\biggl( MTa\frac{\ln n}{n}\frac
{n^{a/n}}{n^{a/n}-1}\biggr).
\end{eqnarray*}
But,
\[
\frac{\ln n}{n}\frac{n^{a/n}}{n^{a/n}-1} \sim\frac{\ln n}{n}\frac
{1}{({a\ln n})/{n}} \sim\frac{1}{a},
\]
when $n \rightarrow+ \infty$. Thus, we have shown the result.

%s6 ###
\section{\texorpdfstring{Proof of Lemma
\protect\lowercase{\ref{lemme technique 2}}}{Proof of Lemma 4.5}}
\label{appB}
Thanks to Lemma \ref{lemme technique 1}, we have
\begin{eqnarray*}
\frac{ \prod_{i=0}^{n-1}(1+M_1h_i+M_2
{h_i}/({T-t_{i+1}}))
}{ \prod_{i=0}^{n-1}(1+M_2{h_i}/({T-t_{i+1}})) } &=&
\prod
_{i=0}^{n-1}\biggl( 1+\frac{M_1}{1+M_2{h_i}/({T-t_{i+1}})}h_i
\biggr)\\
&\leq& \prod_{i=0}^{n-1}(1+M_1h_i) \leq C.
\end{eqnarray*}
So we just have to show that
\[
\prod_{i=0}^{n-1}\biggl(1+M_2\frac{h_i}{T-t_{i+1}}\biggr) \leq Cn^{aM_2}.
\]
But
\[
1+M_2\frac{h_i}{T-t_{i+1}} = 1+M_2(n^{a/n}-1).
\]
So
\begin{eqnarray*}
\prod_{i=0}^{n-1}\biggl(1+M_2\frac{h_i}{T-t_{i+1}}\biggr) &=&
\bigl(1+M_2(n^{a/n}-1)\bigr)^n\\
&=& \exp\biggl( n\ln\biggl(1+aM_2\frac{\ln n}{n}+O\biggl(\frac{\ln^2
n}{n^2}\biggr) \biggr) \biggr)\\
&=& \exp\biggl( aM_2 \ln n + O\biggl(\frac{\ln^2 n}{n}\biggr)
\biggr)
\sim n^{aM_2},
\end{eqnarray*}
when $n \rightarrow+ \infty$. Thus, we have shown the result.

%s7 ###
\section{\texorpdfstring{Proof of Proposition
\protect\lowercase{\ref{u uniformement continue}}}{Proof of Proposition 4.19}}
\label{appC}
We will prove this proposition as the authors of \cite
{Delbaen-Hu-Richou-09} do for their Proposition 4.2. In this proof
we omit the superscript $N,\eps$ for $u$, $Y$ and $Z$ to be more
readable. Let $x_0, x'_0 \in\R^d$ and $t_0, t'_0 \in[0,T]$. By an
argument of symmetry we are allowed to suppose that $t_0 \leq
t'_0$. We have
\[
\vert u(t_0,x_0)-u(t'_0,x'_0)\vert\leq\vert
u(t_0,x_0)-u(t_0,x'_0)\vert+\vert
u(t_0,x'_0)-u(t'_0,x'_0)\vert.
\]
Let us begin with the first term. We will use a classical argument of
linearization:
\begin{eqnarray*}
Y_t^{t_0,x_0}-Y_t^{t_0,x'_0} &=&
g_N(X_T^{t_0,x_0})-g_N(X_T^{t_0,x'_0})\\
&&{}+\int_t^T\alpha_s(X_s^{t_0,x_0}-X_s^{t_0,x'_0})+\beta
_s(Y_s^{t_0,x_0}-Y_s^{t_0,x'_0})\,ds\\
&&{}-\int_t^T (Z_s^{t_0,x_0}-Z_s^{t_0,x'_0})\,d\tilde{W}_s
\end{eqnarray*}
with
\[
\alpha_s:=
\frac{f^{\eps
}(s,X_s^{t_0,x_0},Y_s^{t_0,x'_0},Z_s^{t_0,x'_0})-f^{\eps
}(s,X_s^{t_0,x'_0},Y_s^{t_0,x'_0},Z_s^{t_0,x'_0})}{X_s^{t_0,x_0}-X_s^{t_0,x'_0}},
\]
if $X_s^{t_0,x_0}-X_s^{t_0,x'_0}\neq0$ and $\alpha_s=0$ elsewhere,
\[
\beta_s:=\frac{f^{\eps
}(s,X_s^{t_0,x_0},Y_s^{t_0,x_0},Z_s^{t_0,x'_0})-f^{\eps
}(s,X_s^{t_0,x_0},Y_s^{t_0,x'_0},Z_s^{t_0,x'_0})}{Y_s^{t_0,x_0}-Y_s^{t_0,x'_0}},
\]
if $X_s^{t_0,x_0}-X_s^{t_0,x'_0}\neq0$ and $\beta_s=0$ elsewhere,
\begin{eqnarray*}
\gamma_s&:=&
\frac{f^{\eps}(s,X_s^{t_0,x_0},Y_s^{t_0,x_0},Z_s^{t_0,x_0})-f^{\eps
}(s,X_s^{t_0,x_0},Y_s^{t_0,x_0},Z_s^{t_0,x'_0})}{\vert
Z_s^{t_0,x_0}-Z_s^{t_0,x'_0}\vert^2}\\
&&{}\times\supmathitt{(Z_s^{t_0,x_0}-Z_s^{t_0,x'_0})},
\end{eqnarray*}
if\vspace*{1pt} $Z_s^{t_0,x_0}-Z_s^{t_0,x'_0}\neq0$ and $\gamma_s=0$ elsewhere and
$d\tilde{W}_s:=dW_s-\gamma_s \,ds$.
By a BMO argument, there exists a probability $\Q$ under which $\tilde
{W}$ is a Brownian motion. Then we apply a classical transformation to obtain
\begin{eqnarray*}
&&\E^{\Q}[e^{\int_{t_0}^t\beta_s\,ds}
(Y_t^{t_0,x_0}-Y_t^{t_0,x'_0})] \\
&&\qquad= \E^{\Q}
\biggl[e^{\int
_{t_0}^T\beta_s\,ds}\bigl(g_N(X_T^{t_0,x_0})-g_N(X_T^{t_0,x'_0})
\bigr)
\\
&&\qquad\quad\hspace*{17.5pt}{} + \int_{t_0}^T\alpha_se^{\int_{t_0}^s\beta_u\,du}
(X_s^{t_0,x_0}-X_s^{t_0,x'_0})\,ds\biggr]
\end{eqnarray*}
and
\begin{eqnarray*}
&&\vert u(t_0,x_0)-u(t_0,x'_0)\vert\\
&&\qquad\leq C\biggl( \E^{\Q}[\omega(\vert
X_T^{t_0,x_0}-X_T^{t_0,x'_0}\vert)]+\int_{t_0}^T
\E^{\Q
}
[\vert X_s^{t_0,x_0}-X_s^{t_0,x'_0}\vert]\,ds\biggr)
\end{eqnarray*}
with $\omega$ a modulus of continuity of $g$ that is also a modulus of
continuity for $g_N$. We are allowed to suppose that $\omega$ is
concave; indeed, there exist two positive constants $a$ and $b$ such
that $\omega(x)\leq ax+b$, then the concave hull of $x\mapsto
\omega(x)\vee(\mathbh{1}_{x\geq1}(ax+b))$ is also a
modulus of continuity of $g$. So Jensen's inequality gives us
\begin{eqnarray*}
&&\vert u(t_0,x_0)-u(t_0,x'_0)\vert\\
&&\qquad\leq C\biggl( \omega(\E^{\Q}[\vert
X_T^{t_0,x_0}-X_T^{t_0,x'_0}\vert])+\int_{t_0}^T
\E^{\Q
}
[\vert X_s^{t_0,x_0}-X_s^{t_0,x'_0}\vert]\,ds\biggr).
\end{eqnarray*}
By using the fact that $b$ is bounded we can prove the following
proposition exactly as authors of \cite{Delbaen-Hu-Richou-09} do for
their Proposition 4.7.
\begin{prop}
$\exists C>0$ that does not depend on $N$ and $\eps$ such that
$\forall
t,t' \in[0,T]$, $\forall x,x' \in\R^d$, $\forall s \in[0,T]$,
\[
\E^{\Q}[\vert X_s^{t,x}-X_s^{t',x'}\vert]
\leq C(\vert x-x'\vert+\vert t-t'\vert
^{1/2}).
\]
\end{prop}

Then,
\[
\vert u(t_0,x_0)-u(t_0,x'_0)\vert\leq C\bigl( \omega
(\vert x_0-x'_0\vert)+\vert
x_0-x'_0\vert\bigr).
\]
Now we will study the second term,
\begin{eqnarray*}
\vert u(t_0,x'_0)-u(t'_0,x'_0)\vert&=&\vert
Y_{t_0}^{t_0,x'_0}-Y_{t'_0}^{t'_0,x'_0}\vert\\
& \leq& \vert Y_{t_0}^{t_0,x'_0}-Y_{t_0}^{t'_0,x'_0}
\vert+\vert Y_{t_0}^{t'_0,x'_0}-Y_{t'_0}^{t'_0,x'_0}
\vert.
\end{eqnarray*}
First,
\[
\vert Y_{t_0}^{t'_0,x'_0}-Y_{t'_0}^{t'_0,x'_0}\vert\leq
\biggl\vert\int_{t_0}^{t'_0} f(s,x'_0,Y_s^{t'_0,x'_0},0)\,ds
\biggr\vert\leq C\vert t_0-t'_0\vert.
\]
Moreover, as for the first term we have
\begin{eqnarray*}
&&\E^{\Q}[e^{\int_{t_0}^t\beta_s\,ds}
(Y_t^{t_0,x'_0}-Y_t^{t'_0,x'_0})]\\
&&\qquad= \E^{\Q}\biggl[e^{\int_{t_0}^T\beta_s\,ds}
\bigl(g_N(X_T^{t_0,x'_0})-g_N(X_T^{t'_0,x'_0})\bigr)\\
&&\hspace*{17.1pt}\qquad\quad{}+\int_{t_0}^T\alpha_se^{\int_{t_0}^s\beta_u\,du}
(X_s^{t_0,x'_0}-X_s^{t'_0,x'_0})\,ds\biggr]
\end{eqnarray*}
and
\[
\vert Y_{t_0}^{t_0,x'_0}-Y_{t'_0}^{t'_0,x'_0}\vert\leq
C\bigl( \omega
(\vert t_0-t'_0\vert^{1/2})+\vert
t_0-t'_0\vert^{1/2}\bigr).
\]
Finally,
\[
\vert u(t_0,x'_0)-u(t'_0,x'_0)\vert\leq C\bigl( \omega
(\vert t_0-t'_0\vert^{1/2})+\vert
t_0-t'_0\vert^{1/2}\bigr)
\]
and
\begin{eqnarray*}
&&\vert u(t_0,x_0)-u(t'_0,x'_0)\vert\\
&&\qquad \leq C\bigl(\omega(\vert x_0-x'_0\vert
)+ \omega
(\vert t_0-t'_0\vert^{1/2})+\vert
x_0-x'_0\vert+\vert t_0-t'_0\vert^{1/2}\bigr).
\end{eqnarray*}
So $u$ is uniformly continuous on $[0,T] \times\R^d$ and this function
has a modulus of continuity that does not depend on $N$ and $\eps$.
Moreover, we are allowed to suppose that this modulus of continuity is concave.
\end{appendix}

\section*{Acknowledgments}
The author would like to thank his Ph.D. advisers, Philippe Briand and
Ying Hu, an anonymous referee and the Associate Editor for their
careful reading and helpful comments.

% imsref loaded by lrinkeviciute, 2011-02-19 12:33:37
%
% imsref loaded by lrinkeviciute, 2011-02-19 13:55:29
%

%
\printaddresses


\begin{thebibliography}{22}
% BibTex style file: ims.bst, 2010-01-14
% Default style options (sort=0,type=number).
% Used options (sort=1,type=number).

%b1 ###
\bibitem{Bouchard-Touzi-04}
%
\begin{barticle}[mr]
\bauthor{\bsnm{Bouchard},~\bfnm{Bruno}\binits{B.}} \AND
\bauthor{\bsnm{Touzi},~\bfnm{Nizar}\binits{N.}}
(\byear{2004}).
\btitle{Discrete-time approximation and {M}onte-{C}arlo simulation of backward
stochastic differential equations}.
\bjournal{Stochastic Process. Appl.}
\bvolume{111}
\bpages{175--206}.
\bid{doi={10.1016/j.spa.2004.01.001}, mr={2056536}}
\end{barticle}
%
\endbibitem

%b2 ###
\bibitem{Briand-Confortola-08}
%
\begin{barticle}[mr]
\bauthor{\bsnm{Briand},~\bfnm{Philippe}\binits{P.}} \AND
\bauthor{\bsnm{Confortola},~\bfnm{Fulvia}\binits{F.}}
(\byear{2008}).
\btitle{B{SDE}s with stochastic {L}ipschitz condition and quadratic
{PDE}s in
{H}ilbert spaces}.
\bjournal{Stochastic Process. Appl.}
\bvolume{118}
\bpages{818--838}.
\bid{doi={10.1016/j.spa.2007.06.006}, mr={2411522}}
\end{barticle}
%
\endbibitem

%b3 ###
\bibitem{Briand-Hu-06}
%
\begin{barticle}[mr]
\bauthor{\bsnm{Briand},~\bfnm{Philippe}\binits{P.}} \AND
\bauthor{\bsnm{Hu},~\bfnm{Ying}\binits{Y.}}
(\byear{2006}).
\btitle{B{SDE} with quadratic growth and unbounded terminal value}.
\bjournal{Probab. Theory Related Fields}
\bvolume{136}
\bpages{604--618}.
\bid{doi={10.1007/s00440-006-0497-0}, mr={2257138}}
\end{barticle}
%
\endbibitem

%b4 ###
\bibitem{Briand-Hu-08}
%
\begin{barticle}[mr]
\bauthor{\bsnm{Briand},~\bfnm{Philippe}\binits{P.}} \AND
\bauthor{\bsnm{Hu},~\bfnm{Ying}\binits{Y.}}
(\byear{2008}).
\btitle{Quadratic {BSDE}s with convex generators and unbounded terminal
conditions}.
\bjournal{Probab. Theory Related Fields}
\bvolume{141}
\bpages{543--567}.
\bid{doi={10.1007/s00440-007-0093-y}, mr={2391164}}
\end{barticle}
%
\endbibitem

%b5 ###
\bibitem{Cheridito-Stadje-10}
%
\begin{bmisc}[auto:STB|2010-11-18|09:18:59]
\bauthor{\bsnm{Cheridito},~\bfnm{P.}\binits{P.}} \AND
\bauthor{\bsnm{Stadje},~\bfnm{M.}\binits{M.}}
(\byear{2010}).
\bhowpublished{BS{$\Delta$}Es and BSDEs with non-Lipschitz drivers:
Comparison, convergence and robustness. Available at
\href{http://arxiv.org/abs/1002.0175v1}{arXiv:1002.0175v1}}.
\end{bmisc}
%
\endbibitem

%b6 ###
\bibitem{Chevance-97}
%
\begin{bincollection}[mr]
\bauthor{\bsnm{Chevance},~\bfnm{D.}\binits{D.}}
(\byear{1997}).
\btitle{Numerical methods for backward stochastic differential equations}.
In \bbooktitle{Numerical Methods in Finance}.
\bseries{Publ. Newton Inst.}
\bpages{232--244}.
\bpublisher{Cambridge Univ. Press}, \baddress{Cambridge}.
\bid{mr={1470517}}
\end{bincollection}
%
\endbibitem

%b7 ###
\bibitem{Delarue-Menozzi-06}
%
\begin{barticle}[mr]
\bauthor{\bsnm{Delarue},~\bfnm{Fran{\c{c}}ois}\binits{F.}} \AND
\bauthor{\bsnm{Menozzi},~\bfnm{St{\'e}phane}\binits{S.}}
(\byear{2006}).
\btitle{A forward--backward stochastic algorithm for quasi-linear {PDE}s}.
\bjournal{Ann. Appl. Probab.}
\bvolume{16}
\bpages{140--184}.
\bid{doi={10.1214/105051605000000674}, mr={2209339}}
\end{barticle}
%
\endbibitem

%b8 ###
\bibitem{Delbaen-Hu-Bao-09}
%
\begin{barticle}[auto:STB|2010-11-18|09:18:59]
\bauthor{\bsnm{Delbaen},~\bfnm{F.}\binits{F.}},
\bauthor{\bsnm{Hu},~\bfnm{Y.}\binits{Y.}} \AND
\bauthor{\bsnm{Bao},~\bfnm{X.}\binits{X.}}
(\byear{2011}).
\btitle{Backward SDEs with superquadratic growth}.
\bjournal{Probab. Theory Related Fields}
\bvolume{150}
\bpages{145--192}.
\end{barticle}
%
\endbibitem

%b9 ###
\bibitem{Delbaen-Hu-Richou-09}
%
\begin{bmisc}[auto:STB|2010-11-18|09:18:59]
\bauthor{\bsnm{Delbaen},~\bfnm{F.}\binits{F.}},
\bauthor{\bsnm{Hu},~\bfnm{Y.}\binits{Y.}} \AND
\bauthor{\bsnm{Richou},~\bfnm{A.}\binits{A.}}
(\byear{2011}).
\bhowpublished{On the uniqueness of solutions to quadratic BSDEs with convex
generators and unbounded terminal conditions. \textit{Ann. Inst.
H. Poincar\'e Probab. Statist.} To appear.}
\end{bmisc}
%
\endbibitem

%b10 ###
\bibitem{Fuhrman-Tessitore-02}
%
\begin{barticle}[mr]
\bauthor{\bsnm{Fuhrman},~\bfnm{Marco}\binits{M.}} \AND
\bauthor{\bsnm{Tessitore},~\bfnm{Gianmario}\binits{G.}}
(\byear{2002}).
\btitle{The {B}ismut--{E}lworthy formula for backward {SDE}s and
applications to
nonlinear {K}olmogorov equations and control in infinite dimensional spaces}.
\bjournal{Stochastics Stochastics Rep.}
\bvolume{74}
\bpages{429--464}.
\bid{doi={10.1080/104S1120290024856}, mr={1940495}}
\end{barticle}
%
\endbibitem

%b11 ###
\bibitem{Gobet-Lemor-Warin-05}
%
\begin{barticle}[mr]
\bauthor{\bsnm{Gobet},~\bfnm{Emmanuel}\binits{E.}},
\bauthor{\bsnm{Lemor},~\bfnm{Jean-Philippe}\binits{J.-P.}} \AND
\bauthor{\bsnm{Warin},~\bfnm{Xavier}\binits{X.}}
(\byear{2005}).
\btitle{A regression-based {M}onte {C}arlo method to solve backward stochastic
differential equations}.
\bjournal{Ann. Appl. Probab.}
\bvolume{15}
\bpages{2172--2202}.
\bid{doi={10.1214/105051605000000412}, mr={2152657}}
\end{barticle}
%
\endbibitem

%b12 ###
\bibitem{Gobet-Makhlouf-09}
%
\begin{barticle}[mr]
\bauthor{\bsnm{Gobet},~\bfnm{Emmanuel}\binits{E.}} \AND
\bauthor{\bsnm{Makhlouf},~\bfnm{Azmi}\binits{A.}}
(\byear{2010}).
\btitle{{${\bf L}\sb2$}-time regularity of {BSDE}s with irregular terminal
functions}.
\bjournal{Stochastic Process. Appl.}
\bvolume{120}
\bpages{1105--1132}.
\bid{doi={10.1016/j.spa.2010.03.003}, mr={2639740}}
\end{barticle}
%
\endbibitem

%b13 ###
\bibitem{Hu-Imkeller-Muller-05}
%
\begin{barticle}[mr]
\bauthor{\bsnm{Hu},~\bfnm{Ying}\binits{Y.}},
\bauthor{\bsnm{Imkeller},~\bfnm{Peter}\binits{P.}} \AND
\bauthor{\bsnm{M{\"u}ller},~\bfnm{Matthias}\binits{M.}}
(\byear{2005}).
\btitle{Utility maximization in incomplete markets}.
\bjournal{Ann. Appl. Probab.}
\bvolume{15}
\bpages{1691--1712}.
\bid{doi={10.1214/105051605000000188}, mr={2152241}}
\end{barticle}
%
\endbibitem

%b14 ###
\bibitem{Imkeller-dosReis-09}
%
\begin{barticle}[mr]
\bauthor{\bsnm{Imkeller},~\bfnm{Peter}\binits{P.}} \AND
\bauthor{\bsnm{Dos~Reis},~\bfnm{Gon{\c{c}}alo}\binits{G.}}
(\byear{2010}).
\btitle{Path regularity and explicit convergence rate for {BSDE} with truncated
quadratic growth}.
\bjournal{Stochastic Process. Appl.}
\bvolume{120}
\bpages{348--379}.
\bid{doi={10.1016/j.spa.2009.11.004}, mr={2584898}}
\end{barticle}
%
\endbibitem

%b15 ###
\bibitem{Imkeller-Reis-Zhang-10}
%
\begin{bincollection}[auto:STB|2010-11-18|09:18:59]
\bauthor{\bsnm{Imkeller},~\bfnm{P.}\binits{P.}}, \bauthor{\bparticle{dos}
\bsnm{Reis},~\bfnm{G.}\binits{G.}} \AND
\bauthor{\bsnm{Zhang},~\bfnm{J.}\binits{J.}}
(\byear{2010}).
\btitle{Results on numerics for FBSDE with drivers of quadratic growth}.
In \bbooktitle{Contemporary Quantitative Finance}
(\beditor{A. Chiarella and C.~Novikov}, eds.).
\bpublisher{Springer}, \baddress{Berlin}.
\end{bincollection}
%
\endbibitem

%b16 ###
\bibitem{Kazamaki-94}
%
\begin{bbook}[mr]
\bauthor{\bsnm{Kazamaki},~\bfnm{Norihiko}\binits{N.}}
(\byear{1994}).
\btitle{Continuous Exponential Martingales and {BMO}}.
\bseries{Lecture Notes in Math.}
\bvolume{1579}.
\bpublisher{Springer}, \baddress{Berlin}.
\bid{mr={1299529}}
\end{bbook}
%
\endbibitem

%b17 ###
\bibitem{Kobylanski-00}
%
\begin{barticle}[mr]
\bauthor{\bsnm{Kobylanski},~\bfnm{Magdalena}\binits{M.}}
(\byear{2000}).
\btitle{Backward stochastic differential equations and partial differential
equations with quadratic growth}.
\bjournal{Ann. Probab.}
\bvolume{28}
\bpages{558--602}.
\bid{doi={10.1214/aop/1019160253}, mr={1782267}}
\end{barticle}
%
\endbibitem

%b18 ###
\bibitem{Lepeltier-SanMartin-98}
%
\begin{barticle}[mr]
\bauthor{\bsnm{Lepeltier},~\bfnm{J.~P.}\binits{J.~P.}} \AND
\bauthor{\bsnm{San~Mart{\'{\i}}n},~\bfnm{J.}\binits{J.}}
(\byear{1998}).
\btitle{Existence for {BSDE} with superlinear-quadratic coefficient}.
\bjournal{Stochastics Stochastics Rep.}
\bvolume{63}
\bpages{227--240}.
\bid{mr={1658083}}
\end{barticle}
%
\endbibitem

%b19 ###
\bibitem{Zhang-04}
%
\begin{barticle}[mr]
\bauthor{\bsnm{Zhang},~\bfnm{Jianfeng}\binits{J.}}
(\byear{2004}).
\btitle{A numerical scheme for {BSDE}s}.
\bjournal{Ann. Appl. Probab.}
\bvolume{14}
\bpages{459--488}.
\bid{doi={10.1214/aoap/1075828058}, mr={2023027}}
\end{barticle}
%
\endbibitem

%b20 ###
\bibitem{Zhang-05}
%
\begin{barticle}[mr]
\bauthor{\bsnm{Zhang},~\bfnm{Jianfeng}\binits{J.}}
(\byear{2005}).
\btitle{Representation of solutions to {BSDE}s associated with a degenerate
{FSDE}}.
\bjournal{Ann. Appl. Probab.}
\bvolume{15}
\bpages{1798--1831}.
\bid{doi={10.1214/105051605000000232}, mr={2152246}}
\end{barticle}
%
\endbibitem

%b21 ###
\bibitem{Zhang-01}
%
\begin{bmisc}[auto:STB|2010-11-18|09:18:59]
\bauthor{\bsnm{Zhang},~\bfnm{J.}\binits{J.}}
(\byear{2001}).
\bhowpublished{Some fine properties of BSDE. Ph.D. thesis, Purdue
Univ.}
\end{bmisc}
%
\endbibitem

\end{thebibliography}
\end{document}